\documentclass[a4paper,12pt]{article}

\usepackage{amssymb}                                    
\usepackage{mathrsfs}                    
\usepackage{amsmath}                    
\usepackage{amsfonts}                   
\usepackage{amsthm}                     
\usepackage[latin1]{inputenc}           
\usepackage[arrow, matrix, curve]{xy}    
\usepackage[english]{babel}                
\usepackage{graphics}
\usepackage{color}
\usepackage{hyperref}
\usepackage[normalem]{ulem}

\usepackage{shorttoc}

\theoremstyle{plain}
\newtheorem{theorem}{Theorem}[section]
\newtheorem{corollary}[theorem]{Corollary}
\newtheorem{proposition}[theorem]{Proposition}
\newtheorem{lemma}[theorem]{Lemma}

\newtheorem{observation}[theorem]{Observation}

\theoremstyle{definition}
\newtheorem{definition}[theorem]{Definition}
\newtheorem{remark}[theorem]{Remark}
\newtheorem{example}[theorem]{Example}

\newcommand{\tfleft}{\!\xymatrix@C=1.3pc{&\ar@{->>}_\sim[l]}\!}

\newcommand{\tfright}{\!\xymatrix@C=1.3pc{\ar@{->>}^\sim[r]&}\!}

\def\endofproof{\hfill{$\square$}\\}

\title{Principal $\infty$-bundles -- General theory }
\author{Thomas Nikolaus, Urs Schreiber, Danny Stevenson}
\date{\today}

\begin{document}

\maketitle

\begin{abstract}
 The theory of principal bundles
 makes sense in any $\infty$-topos, such as the $\infty$-topos
 of topological, of smooth, or of otherwise geometric $\infty$-groupoids/$\infty$-stacks,
 and more generally in slices of these.
 It provides a natural
 geometric model for structured higher nonabelian cohomology and controls general
 fiber bundles in terms of associated bundles.
 For suitable choices of structure $\infty$-group $G$ these \emph{$G$-principal $\infty$-bundles}
 reproduce various higher structures that have been considered in the literature and
 further generalize these to a full geometric model for twisted higher nonabelian sheaf cohomology.
 We discuss here this general abstract theory of principal $\infty$-bundles, observing that
 it is intimately related to the axioms
 that characterize $\infty$-toposes.
 A central result is a natural
 equivalence between principal $\infty$-bundles and intrinsic nonabelian cocycles,
 implying the classification of principal $\infty$-bundles by
 nonabelian sheaf hyper-cohomology. We observe that the theory of geometric fiber $\infty$-bundles
 \emph{associated} to principal $\infty$-bundles subsumes a theory of \emph{$\infty$-gerbes}
 and of \emph{twisted $\infty$-bundles}, with twists deriving
 from \emph{local coefficient $\infty$-bundles},  which we define,
 relate to extensions of principal $\infty$-bundles and show to be classified by a
 corresponding notion of \emph{twisted cohomology}, identified with the
 cohomology of a corresponding slice $\infty$-topos.


 \end{abstract}

\newpage

\tableofcontents

\newpage

\section{Overview}
\label{Overview}

The concept of a $G$-\emph{principal bundle} for a topological or
Lie group $G$ is fundamental in
classical topology and differential
geometry, e.g. \cite{Husemoeller}. More generally, for $G$ a \emph{geometric}
group in the sense of a \emph{sheaf of groups} over some site, the notion of
$G$-principal bundle or \emph{$G$-torsor} is fundamental in \emph{topos theory}
\cite{Johnstone,Moerdijk}.
Its relevance rests in the fact that
$G$-principal bundles constitute natural geometric representatives
of cocycles in degree 1 nonabelian cohomology
$H^1(-,G)$  and that general fiber bundles are
\emph{associated} to principal bundles.

In recent years it has become clear that
various applications, notably in
``\emph{String-geometry}'' \cite{SSS,Schreiber},  involve a notion of principal bundles where
geometric groups $G$ are generalized to \emph{geometric grouplike $A_\infty$-spaces},
in other words
\emph{geometric $\infty$-groups}: geometric objects that are equipped with a group
structure up to higher coherent homotopy. The resulting \emph{principal $\infty$-bundles}
should be natural geometric representatives of geometric nonabelian \emph{hypercohomology}:
{\v C}ech cohomology with coefficients in arbitrary positive degree.

In the \emph{absence} of geometry,
these principal $\infty$-bundles are essentially just the classical \emph{simplicial} principal bundles
of simplicial sets \cite{May} (this we discuss in Section 4.1 of \cite{NSSb}). However,
in the presence of non-trivial geometry the situation is both more subtle and richer, and plain
simplicial principal bundles can only serve as a specific \emph{presentation}
for the general notion (section 3.7.2 of \cite{NSSb}).

For the case of \emph{principal 2-bundles},
which is the first step after ordinary principal bundles,
aspects of a geometric definition and theory have been proposed and developed by
various authors, see section 1 of \cite{NSSb} for references and see
\cite{NikolausWaldorf} for a comprehensive discussion.
Notably the notion of a \emph{bundle gerbe} \cite{Mur}
is, when regarded as an extension of a {\v C}ech-groupoid,
almost manifestly that of a principal 2-bundle,
even though this perspective is not prominent in the respective literature.
The oldest definition of geometric 2-bundles is conceptually different, but closely related:
Giraud's \emph{$G$-gerbes} \cite{Giraud} are by definition not principal 2-bundles but are
fiber 2-bundles \emph{associated} to $\mathbf{Aut}(\mathbf{B}G)$-principal 2-bundles,
where $\mathbf{B}G$ is the \emph{geometric moduli stack} of $G$-principal bundles.
This means that $G$-gerbes provide the universal \emph{local coefficients}, in the sense
of \emph{twisted cohomology}, for $G$-principal bundles.

From the definition of principal 2-bundles/bundle gerbes it is fairly clear that these
ought to be just the first step (or second step) in an infinite tower of higher analogs.
Accordingly, definitions of \emph{principal 3-bundles} have
been considered in the literature, mostly in the guise of
\emph{bundle 2-gerbes} \cite{Stevenson}
.
The older notion
of Breen's \emph{$G$-2-gerbes} \cite{Breen} (also discussed by Brylinski-MacLaughlin),
is, as before, not that of a principal
3-bundle, but that of a fiber 3-bundle which is
\emph{associated} to an $\mathbf{Aut}(\mathbf{B}G)$-principal 3-bundle,
where now
$\mathbf{B}G$ is the \emph{geometric moduli 2-stack} of $G$-principal 2-bundles
.

Generally, for every $n \in \mathbb{N}$ and every geometric $n$-group $G$,
it is natural to consider the theory of
$G$-principal $n$-bundles \emph{twisted} by an $\mathbf{Aut}(\mathbf{B}G)$-principal
$(n+1)$-bundle, hence by the associated \emph{$G$-$n$-gerbe}.
A complete theory of principal bundles therefore needs to involve
the notion of principal $n$-bundles
and also that of twisted principal $n$-bundles
in the limit as $n \to \infty$.

As $n$ increases, the  piecemeal
conceptualization of principal $n$-bundles quickly becomes tedious
and their structure opaque, without a general theory of higher geometric
structures. In recent years such a theory -- long conjectured and with many precursors --
has materialized in a comprehensive and elegant form, now known as \emph{$\infty$-topos theory}
\cite{ToenVezzosiToposb, Rezk, Lurie}.
Whereas an ordinary topos is a category of \emph{sheaves} over
some site\footnote{Throughout \emph{topos} here stands for \emph{Grothendieck topos},
as opposed to the more general notion of \emph{elementary topos}.}, an
$\infty$-topos is an \emph{$\infty$-category} of \emph{$\infty$-sheaves} or equivalently
of \emph{$\infty$-stacks} over some \emph{$\infty$-site}, where the prefix ``$\infty-$''
indicates that all these notions are generalized to structures up to
\emph{coherent higher homotopy} (as in the older terminology of
$A_\infty$-, $C_\infty$-, $E_\infty$- and $L_\infty$-algebras, all of which re-appear
as algebraic structures in $\infty$-topos theory).
In as far as an ordinary topos is a context
for general \emph{geometry}, an $\infty$-topos is a context for what is called
\emph{higher geometry} or \emph{derived geometry}: the pairing of the notion of
\emph{geometry} with that of \emph{homotopy}.
(Here ``derived'' alludes to
``derived category'' and ``derived functor'' in homological algebra,
but refers in fact to a nonabelian generalization of these concepts.)
Therefore we may refer to objects of an $\infty$-topos also as
\emph{geometric homotopy types}.

As a simple instance of this pairing,
one observes that
for any geometric abelian group (sheaf of abelian groups)
$A$, the higher degree (sheaf) cohomology $H^{n+1}(-,A)$ in ordinary geometry
may equivalently be understood as the degree-1 cohomology
$H^1(-,\mathbf{B}^{n}A)$ in higher geometry, where
$\mathbf{B}^{n}A$ is the geometric $\infty$-group obtained by
successively \emph{delooping} $A$ geometrically.
More generally, there are geometric $\infty$-groups $G$ not of this abelian form.
The general degree-1 geometric cohomology
$H^1(X,G)$ is a nonabelian and simplicial generalization of \emph{sheaf hypercohomology},
whose cocycles are morphisms $X \to \mathbf{B}G$ into the geometric delooping of $G$.
Indeed, delooping plays a central role
in $\infty$-topos theory;
a fundamental fact of $\infty$-topos theory (recalled as
Theorem \ref{DeloopingTheorem} below) says that, quite generally, under internal
looping and delooping, $\infty$-groups $G$ in an $\infty$-topos $\mathbf{H}$ are equivalent
to connected and \emph{pointed} objects in $\mathbf{H}$:
$$
  \xymatrix{
    \left\{
	  \mbox{
	  groups in $\mathbf{H}$
	  }
	\right\}
		\ar@{<-}@<+5pt>[rrr]^<<<<<<<<<<<<<{\mbox{\small  looping}\; \Omega}
	\ar@<-5pt>[rrr]_<<<<<<<<<<<<<{\mbox{\small delooping}\; \mathbf{B}}^<<<<<<<<<<<<<\simeq
    &&&
	\left\{
	  \mbox{\begin{tabular}{c} pointed connected \\ objects in $\mathbf{H}$ \end{tabular}}
	\right\}
  }
  \,.
$$
We will see that this
equivalence of $\infty$-categories
plays a key role in the theory of
principal $\infty$-bundles.

Topos theory is renowned for providing a general convenient context for the development
of geometric structures. In some sense, $\infty$-topos theory
provides an even more convenient context, due to the fact that
\emph{$\infty$-(co)limits} or \emph{homotopy (co)limits} in an $\infty$-topos
exist, and refine
the corresponding naive (co)limits.
This convenience manifests itself in the central definition of
principal $\infty$-bundles (Definition~\ref{principalbundle} below): whereas the
traditional definition of a $G$-principal bundle over $X$ as a quotient map
$P \to P/G \simeq X$ requires the additional clause that the quotient be
\emph{locally trivial}, $\infty$-topos theory comes pre-installed with the correct homotopy quotient for
higher geometry, and as a result the local triviality of $P \to P/\!/G =: X$ is automatic;
we discuss this in more detail in Section \ref{PrincBund-Intro} below.
Hence
conceptually, $G$-principal $\infty$-bundles are in fact simpler than their
traditional analogs, and so their theory is stronger.

A central theorem of topos theory is \emph{Giraud's theorem}, which
intrinsically characterizes
toposes as those presentable categories that satisfy three simple
conditions: 1.\  coproducts are disjoint, 2.\  colimits are preserved by pullback, and
3.\  quotients are effective. The analog of this characterization turns out
to remain true essentially verbatim in $\infty$-topos theory:
this is the \emph{Giraud-To{\"e}n-Vezzosi-Rezk-Lurie} characterization
of $\infty$-toposes, recalled as Definition \ref{GiraudRezkLurieAxioms} below.
We will show that given an $\infty$-topos $\mathbf{H}$, the second
and the third of these axioms
lead directly to the {\em classification theorem} for principal $\infty$-bundles
(Theorem~\ref{PrincipalInfinityBundleClassification} below) which states that
there is an equivalence of $\infty$-groupoids
\[
\begin{array}{cccc}
G \mathrm{Bund}(X)
	  &
	  \simeq
	  &
	  \mathbf{H}(X, \mathbf{B}G)
\\
\\
  	\left\{
	  \begin{array}{c}
	  \mbox{$G$-principal $\infty$-bundles}
	  \\
	  \mbox{over $X$}
	  \end{array}
	  \raisebox{20pt}{
	  \xymatrix{
	     P \ar[d] \\ X
	  }}
	\right\}
    &\simeq&
  	\left\{
	  \mbox{cocycles}\;\;\;
	   g : X \to \mathbf{B}G
	\right\}
	  	\end{array}
\]	
between the $\infty$-groupoid of $G$-principal $\infty$-bundles on $X$,
and the mapping space $\mathbf{H}(X,\mathbf{B}G)$.

The mechanism underlying the proof of this theorem is summarized
in the following diagram, which is supposed to indicate that the geometric $G$-principal $\infty$-bundle corresponding
to a cocycle is nothing but the
corresponding homotopy fiber:
 $$
  \raisebox{50pt}{
  \xymatrix@C=-7pt@R=2pt{
    \vdots && \vdots
    \\
    P \times G \times G
	  \ar[rr] \ar@<-4pt>[ddd] \ar@<+0pt>[ddd] \ar@<+4pt>[ddd]
	  &&
	G \times G
      \ar@<-4pt>[ddd] \ar@<+0pt>[ddd] \ar@<+4pt>[ddd]	
    \\
	\\
    \\
    P \times G \ar[rr] \ar@<-3pt>[ddd]_{p_1} \ar@<+3pt>[ddd]^\rho
	&& G
	\ar@<-3pt>[ddd]_{} \ar@<+3pt>[ddd]
    \\
	\\ && & \mbox{\small $G$-$\infty$-actions}
	\\
    P \ar[rr] \ar[ddd] && {*} \ar[ddd] & \mbox{\small total objects}
	\\
	\\ & \mbox{\small $\infty$-pullback} &&
	\\
	X \ar[rr]_{g } &&  \mathbf{B}G &  \mbox{\small quotient objects}
	\\
	\mbox{\small \begin{tabular}{c} $G$-principal \\ $\infty$-bundle \end{tabular}}
	\ar@{}[rr]_{\mbox{\small cocycle}}
	&& \mbox{\small \begin{tabular}{c} universal \\ $\infty$-bundle \end{tabular}}
  }
  }
$$
The fact that all geometric $G$-principal $\infty$-bundles arise this
way, up to equivalence, is quite useful in applications, and also sheds helpful light
on various existing constructions and provides more examples.

Notably, the implication
that \emph{every} geometric $\infty$-action $\rho : V \times G \to V$
of an $\infty$-group $G$ on an object $V$
has a classifying morphism $\mathbf{c} : V/\!/G \to \mathbf{B}G$,
tightly connects the theory of associated $\infty$-bundles with that of
principal $\infty$-bundles (Section~\ref{StrucRepresentations} below):
the fiber sequence
$$
  \raisebox{20pt}{
  \xymatrix{
    V \ar[r] & V/\!/G
	\ar[d]^{\mathbf{c}}
	\\
	& \mathbf{B}G
  }
  }
$$
is found to be the $V$-fiber $\infty$-bundle which is \emph{$\rho$-associated}
to the universal $G$-principal $\infty$-bundle $* \to \mathbf{B}G$.
Again, using the $\infty$-Giraud axioms,
an $\infty$-pullback of $\mathbf{c}$ along a cocycle $g_X : X \to \mathbf{B}G$
is identified with the $\infty$-bundle $P \times_G V$ that is \emph{$\rho$-associated}
to the principal $\infty$-bundle $P \to X$ classified by $g_X$
(Proposition \ref{RhoAssociatedToUniversalIsUniversalVBundle}) and
every $V$-fiber $\infty$-bundle arises this way, associated to an
$\mathbf{Aut}(V)$-principal $\infty$-bundle (Theorem \ref{VBundleClassification}).

Using this, we may observe that the space $\Gamma_X(P \times_G V)$
of \emph{sections} of $P \times_G V$ is equivalently
the space $\mathbf{H}_{/\mathbf{B}G}(g_X, \mathbf{c})$
of cocycles $\sigma : g_X \to \mathbf{c}$ in the slice $\infty$-topos
$\mathbf{H}_{/\mathbf{B}G}$:
$$
 \begin{array}{cccc}
  \Gamma_X(P \times_G V)
	  &
	  \simeq
	  &
	  \mathbf{H}_{/\mathbf{B}G}(g_X, \mathbf{c})
	   \\
	  \\
  	\left\{
	  \raisebox{20pt}{
	  \xymatrix{
	    P \times_G V \ar[r] \ar[d] & V /\! / G
		  \ar[d]^-{\mathbf{c}}
	    \\
	    X \ar[r]_-{g_X}
		\ar@/^1pc/@{-->}[u]^{\sigma}
		& \mathbf{B}G
	  }
	  }
	\right\}
    &\simeq&
  	\left\{
	  \raisebox{20pt}{
	  \xymatrix{
	    & V /\! / G
		  \ar[d]^-{\mathbf{c}}
	    \\
	    X \ar[r]_-{g_X}
		\ar@{-->}[ur]^{\sigma}
		& \mathbf{B}G
	  }
	  }
	\right\}
	\end{array}
$$
Moreover, by the above classification theorem of $G$-principal $\infty$-bundles,
$g_X$ trivializes over some
cover $\xymatrix{U \ar@{->>}[r] &  X}$,
and so the universal property of the $\infty$-pullback implies that
\emph{locally} a section $\sigma$ is a $V$-valued function
$$
 \raisebox{20pt}{
  \xymatrix{
    & V \ar[r] & V /\! / G
	  \ar[d]^{\mathbf{c}}
    \\
    U \ar@{->>}[r]_<<<{\mathrm{cover}}\ar@{-->}[ur]^{\sigma|_U} & X \ar[r]^{g_X}
	& \mathbf{B}G\, .
  }
  }
$$
For $V$ an ordinary space, hence a \emph{0-truncated} object in the $\infty$-topos,
this is simply the familiar statement about sections of associated bundles.
But in higher geometry $V$,
being an object of the $\infty$-topos and hence an $\infty$-stack,
may more generally itself be a higher
moduli $\infty$-stack,
classifying some geometric structures,
which makes the general theory of sections more interesting.
Specifically, if $V$ is a pointed connected object, then it is of the form
$\mathbf{B}\Omega V$ and  this
means that it is locally a cocycle for an $\Omega V$-principal $\infty$-bundle,
and so globally is a \emph{twisted $\Omega V$-principal $\infty$-bundle}.
This identifies $\mathbf{H}_{/\mathbf{B}G}(-, \mathbf{c})$ as the
\emph{twisted cohomology} induced by the \emph{local coefficient bundle $\mathbf{c}$}
with \emph{local coefficients} $V$.
This yields a geometric and unstable analogue
of the
picture of twisted cohomology discussed in \cite{AndoBlumbergGepner}.

Given $V$, the most general
twisting group is the \emph{automorphism $\infty$-group}
$\mathbf{Aut}(V) \hookrightarrow [V,V]_{\mathbf{H}}$,
formed in the $\infty$-topos (Definition \ref{InternalAutomorphismGroup}). If
$V$ is pointed connected and hence of the form $V = \mathbf{B}G$, this means that the
most general universal local coefficient bundle is
$$
  \raisebox{20pt}{
  \xymatrix{
    \mathbf{B}G \ar[r]
	  &
	 (\mathbf{B}G)/\!/\mathbf{Aut}(\mathbf{B}G)
	 \ar[d]^{\mathbf{c}_{\mathbf{B}G}}
	 \\
	  & \mathbf{B}\mathbf{Aut}(\mathbf{B}G).
  }}
$$
The corresponding associated twisting $\infty$-bundles are
\emph{$G$-$\infty$-gerbes}: fiber $\infty$-bundles with typical fiber
the moduli $\infty$-stack $\mathbf{B}G$. These are the
universal local coefficients for twists of $G$-principal $\infty$-bundles.

While twisted cohomology in $\mathbf{H}$ is hence identified simply with ordinary cohomology
in a slice of $\mathbf{H}$, the corresponding geometric representatives, the
$\infty$-bundles, do not translate to the
slice quite as directly. The reason is that a universal local coefficient bundle
$\mathbf{c}$ as above is rarely a pointed connected object in the slice
(if it is, then it is essentially trivial) and so the theory of principal
$\infty$-bundles does not directly apply to these coefficients.
In Section~\ref{ExtensionsOfCohesiveInfinityGroups} we show that what
does translate is a notion of
\emph{twisted $\infty$-bundles}, a generalization of the twisted bundles known from
twisted K-theory: given a section $\sigma : g_X \to \mathbf{c}$ as above,
the following pasting diagram of $\infty$-pullbacks
$$
  \xymatrix{
    Q
	\ar[d] \ar[r] & {*} \ar[d]
	&&
	\mbox{$P$-twisted $\Omega V$-principal $\infty$-bundle}
    \\
    P \ar[d] \ar[r] & V \ar[r]  \ar[d] & {*} \ar[d] & \mbox{$G$-principal $\infty$-bundle}
    \\
    X \ar[r]^\sigma \ar@/_2pc/[rr]_{g_X} & V/\!/G
    \ar[r]^{\mathbf{c}} & \mathbf{B}G & \mbox{section of $\rho$-associated $V$-$\infty$-bundle}
  }
$$
naturally identifies an $\Omega V$-principal $\infty$-bundle $Q$ on the total space $P$
of the twisting $G$-principal $\infty$-bundle, and since this is classified by
a $G$-equivariant morphism $P \to V$ it enjoys itself a certain twisted $G$-equivariance
with respect to the defining $G$-action on $P$. We call such $Q \to P$ the
\emph{$[g_X]$-twisted $\Omega V$-principal bundle} classified by $\sigma$.
Again, a special case of special importance is that where
$V = \mathbf{B}A$ is pointed connected,
which identifies the universal $V$-coefficient bundle with an
\emph{extension of $\infty$-groups}
$$
  \raisebox{20pt}{
  \xymatrix{
    \mathbf{B}A \ar[r] & \mathbf{B}\hat G
	\ar[d]
	\\
	& \mathbf{B}G.
  }
  }
$$
Accordingly, $P$-twisted $A$-principal $\infty$-bundles are equivalently
\emph{extensions} of $P$ to $\hat G$-principal $\infty$-bundles.

A direct generalization of the previous theorem yields the
classification Theorem \ref{ClassificationOfTwistedGEquivariantBundles},
which identifies $[g_X]$-twisted $A$-principal $\infty$-bundles with cocycles in twisted cohomology
$$
 \begin{array}{cccc}
 	  A \mathrm{Bund}^{[g_X]}(X)
	  &
	  \simeq
	  &
	  \mathbf{H}_{/\mathbf{B}G}(g_X, \mathbf{c}_{\mathbf{B}G})
      \\
      \\
  	\left\{
	  \begin{array}{c}
	   \mbox{$[g_X]$-twisted}
	   \\
	  \mbox{$A$-principal $\infty$-bundles}
	  \\
	  \mbox{over $X$}
	  \end{array}
	  \raisebox{40pt}{
	  \xymatrix{
	     Q \ar[d] \\ P = (g_X)^* {*} \ar[d] \\ X
	  }}
	\right\}
    &\simeq&
  	\left\{
	  \mbox{twisted cocycles}\;\;\;
	   \sigma : g_X \to \mathbf{c}_{\mathbf{B}G}
	\right\}
	\end{array}
$$
For instance if $\mathbf{c}$ is the connecting homomorphism
$$
  \xymatrix{
    \mathbf{B}\hat G \ar[r] & \mathbf{B}G \ar[d]^{\mathbf{c}}
	\\
	& \mathbf{B}^2 A
  }
$$
of a central
extension of ordinary groups $A \to \hat G \to G$, then
the corresponding twisted $\hat G$-bundles are those known from geometric models of
twisted K-theory.

When the internal \emph{Postnikov tower} of a coefficient object is regarded as
a sequence of local coefficient bundles as above, the induced
twisted $\infty$-bundles are decompositions of nonabelian principal $\infty$-bundles
into ordinary principal bundles together with
equivariant abelian hypercohomology cocycles on their total spaces.
This construction identifies much of equivariant cohomology theory
as a special case of higher nonabelian cohomology.
Specifically, when applied to a Postnikov stage of the delooping of an
$\infty$-group of internal automorphisms, the corresponding
twisted cohomology reproduces
the notion of Breen \emph{$G$-gerbes with band} (Giraud's \emph{lien}s);
and the corresponding
twisted $\infty$-bundles are their incarnation as
equivariant \emph{bundle gerbes} over principal bundles.

\medskip

The classification statements for principal and fiber $\infty$-bundles in this article,
Theorems \ref{PrincipalInfinityBundleClassification}
and \ref{VBundleClassification} are not surprising, they say exactly what one would
hope for. It is however useful to see how they flow naturally from the abstract
axioms of $\infty$-topos theory, and to observe that they immediately imply
a series of classical as well as recent theorems as special cases,
see Remark \ref{ReferencesOnClassificationOfVBundles}.
Also the corresponding long exact sequences in (nonabelian) cohomology,
Theorem \ref{LongExactSequenceInCohomology}, reproduce classical theorems,
see Remark \ref{ReferencesOnLongSequences}.
Similarly the
definition and classification of lifting of principal $\infty$-bundles,
Theorem \ref{ExtensionsAndTwistedCohomology}, and of twisted principal $\infty$-bundles
in Theorem \ref{ClassificationOfTwistedGEquivariantBundles} flows naturally
from the $\infty$-topos theory and yet it immediately implies various
constructions and results in the literature as special cases, see
Remark \ref{ReferencesOnLiftings} and Remark \ref{ReferencesTwistedBundles}, respectively.
In particular the notion of nonabelian twisted cohomology itself is elementary
in $\infty$-topos theory,
Section \ref{TwistedCohomology}, and yet it sheds light on  a wealth of
applications, see Remark \ref{ReferencesOnTwisted}.

This should serve to indicate that the theory of (twisted) principal $\infty$-bundles
is rich and interesting.
The present article is intentionally written in general abstraction only, aiming to present
the general theory of (twisted) principal $\infty$-bundles as elegantly as possible,
true to its roots in
abstract higher topos theory. We believe that this serves to usefully make transparent
the overall picture. In the companion article
\cite{NSSb} we give a complementary discussion and construct
explicit presentations of the structures appearing here that lend themselves
to explicit computations.

\section{Preliminaries}
\label{structures}

The discussion of principal $\infty$-bundles in
Section~\ref{Principal infinity-bundles general abstract}
below builds on
the concept of an \emph{$\infty$-topos} and on a handful
of basic structures and notions that are present in any $\infty$-topos,
in particular the notion of \emph{group objects} and of \emph{cohomology} with coefficients
in these group objects. The relevant theory has been
developed in \cite{ToenVezzosiToposb, Rezk, Lurie, LurieAlgebra}.
The purpose of this section
is to recall the main aspects of this theory that we need, to establish our notation, and to
highlight some aspects of the general theory
that are relevant to our discussion and which have perhaps
not been highlighted in this way in the
existing literature.

\medskip

One may reason about $\infty$-categories in a model-independent way, using the universal
constructions that hold equivalently in all models -- but the reader wishing to do so is invited to think
specifically of the model given by quasi-categories, due to Joyal and laid out in detail in
\cite{Lurie}. Ordinary categories naturally embed into $\infty$-categories; and in
terms of quasi-categories this embedding is given by sending a category to its
simplicial nerve. In view of this we will freely regard 1-categories as $\infty$-categories --
such as for instance the simplex category $\Delta$. This allows
us to define a simplicial object in an $\infty$-category $\mathcal{C}$ in a direct generalization of the usual
notion as an $\infty$-functor $\Delta^{\mathrm{op}} \to \mathcal{C}$ (note that
this is now automatically a simplicial object ``up to coherent higher homotopy'').

We are particularly concerned with those $\infty$-categories that are
$\infty$-toposes.
For many purposes the notion of \emph{$\infty$-topos} is best thought of as
a generalization of the notion of a sheaf topos --- the category
of sheaves over some site is replaced by an
\emph{$\infty$-category of $\infty$-stacks/$\infty$-sheaves}
over some $\infty$-site (there is also supposed to be a more general notion of an
\emph{elementary $\infty$-topos}, which however we do not consider here).
In this context the $\infty$-topos $\mathrm{Gpd}_\infty$ of $\infty$-groupoids
is the natural generalization of the punctual topos $\mathrm{Set}$ to
$\infty$-topos theory.  A major achievement
of \cite{ToenVezzosiToposb}, \cite{Rezk} and \cite{Lurie} was to provide a more intrinsic characterization of
$\infty$-toposes, which generalizes the classical characterization of
sheaf toposes (Grothendieck toposes) originally given by Giraud.
We will show that the theory of
principal $\infty$-bundles is naturally expressed in terms of these intrinsic
properties, and therefore we here take these \emph{Giraud-To{\"e}n-Vezzosi-Rezk-Lurie axioms}
to be the very definition of an $\infty$-topos
(\cite{Lurie}, Theorem 6.1.0.6, the main ingredients will be recalled below):
\begin{definition}[$\infty$-Giraud axioms]
  \label{GiraudRezkLurieAxioms}
  An \emph{$\infty$-topos} is a presentable $\infty$-category $\mathbf{H}$ that satisfies the
  following properties.
  \begin{enumerate}
    \item {\bf Coproducts are disjoint.}  For every two objects
    $A, B \in \mathbf{H}$, the intersection of $A$ and
    $B$ in their coproduct is the initial object: in other words the diagram
	$$
	  \xymatrix{
	    \emptyset \ar[r] \ar[d] & B \ar[d]
		\\
		A \ar[r] & A \coprod B
	  }
	$$
	is a pullback.
	
	\item {\bf Colimits are preserved by pullback.}
	  For all morphisms $f\colon X\to B$ in $\mathbf{H}$ and all
	  small diagrams $A\colon I\to \mathbf{H}_{/B}$, there is an
	  equivalence
	  $$
	    \varinjlim_i f^*A_i \simeq f^*(\varinjlim_i A_i)
	  $$	
	  between the pullback of the colimit and the colimit over the pullbacks of its
	  components.
	\item
	  {\bf Quotient maps are effective epimorphisms.} Every simplicial object
	  $A_\bullet : \Delta^{\mathrm{op}} \to \mathbf{H}$ that satisfies the
	  groupoidal Segal property (Definition~\ref{GroupoidObject}) is the {\v C}ech nerve of its quotient projection:
	  $$
	    A_n \simeq
		A_0 \times_{\varinjlim_n A_n} A_0 \times_{\varinjlim_n A_n} \cdots \times_{\varinjlim_n A_n} A_0
		\;\;\;
		\mbox{($n$ factors)}
		\,.
	  $$
  \end{enumerate}
\end{definition}
Repeated application of the second and third axiom provides the proof of
the classification of principal $\infty$-bundles,
Theorem \ref{PrincipalInfinityBundleClassification} and the
universality of the universal associated $\infty$-bundle,
Proposition \ref{UniversalAssociatedBundle}.

An ordinary topos is famously characterized by the existence of a classifier object
for monomorphisms, the \emph{subobject classifier}.
With hindsight, this statement already carries in it the seed
of the close relation between topos theory and bundle theory, for we may think of
a monomorphism $E \hookrightarrow X$ as being a \emph{bundle of $(-1)$-truncated fibers}
over $X$. The following axiomatizes the existence of arbitrary universal bundles,
providing a different but equivalent definition of $\infty$-toposes.
\begin{proposition}
  An \emph{$\infty$-topos} $\mathbf{H}$ is a presentable $\infty$-category with the following
  properties.
  \begin{enumerate}
	\item {\bf Colimits are preserved by pullback.}
	\item {\bf There are universal $\kappa$-small bundles.}
	  For every sufficiently large regular cardinal $\kappa$, there exists
	  a morphism $\widehat {\mathrm{Obj}}_\kappa \to \mathrm{Obj}_\kappa$ in $\mathbf{H}$,
	  such that for every other object $X$, pullback along morphisms $X \to \mathrm{Obj}$
	  constitutes an equivalence\footnote{Here $\mathrm{Core}$ denotes the maximal $\infty$-groupoid
inside an $\infty$-category.}
	  $$
	    \mathrm{Core}(\mathbf{H}_{/_{\kappa}X})
		\simeq
		\mathbf{H}(X, \mathrm{Obj}_\kappa)
	  $$
	  between the $\infty$-groupoid of bundles (morphisms) $E \to X$ which are $\kappa$-small over $X$
	  and the $\infty$-groupoid of morphisms from $X$ into $\mathrm{Obj}_\kappa$.	
  \end{enumerate}
  \label{RezkCharacterization}
\end{proposition}
This is due to Rezk and Lurie, appearing as Theorem 6.1.6.8 in \cite{Lurie}.
We find that this second version of the axioms
naturally gives the equivalence between $V$-fiber bundles and
$\mathbf{Aut}(V)$-principal $\infty$-bundles
in Proposition \ref{VBundleIsAssociated}.

In addition to these axioms, a basic property of $\infty$-toposes
(and generally of $\infty$-categories with pullbacks) which we will
repeatedly invoke, is the following.
\begin{proposition}[pasting law for pullbacks]
 Let $\mathbf{H}$ be an $\infty$-category with pullbacks.  If
  $$
    \xymatrix{
	  A \ar[r]\ar[d] & B \ar[r]\ar[d] & C \ar[d]
	  \\
	  D \ar[r] & E \ar[r] & F
	}
  $$
  is a diagram in $\mathbf{H}$ such that the right
  square is an $\infty$-pullback, then the left square is an $\infty$-pullback precisely
  if the outer rectangle is.
  \label{PastingLawForPullbacks}
\end{proposition}
Notice that here and in all of the following
\begin{itemize}
  \item all square diagrams are filled by a 2-cell, even if we do not indicate this
    notationally;
  \item
    all limits are $\infty$-limits/homotopy limits
	(hence all pullbacks are $\infty$-pullbacks/homotopy pullbacks), and so on;
\end{itemize}
this is the only consistent way of speaking about $\mathbf{H}$ in generality.
Only in the followup article \cite{NSSb} do we consider presentations of
$\mathbf{H}$ by 1-categorical data; there we will draw a careful
distinction between 1-categorical
limits and $\infty$-categorical/homotopy limits.
\begin{definition}
  For $f : Y \to Z$ any morphism in $\mathbf{H}$
  and $z : * \to Z$ a point, the \emph{fiber}
  (homotopy fiber or $\infty$-fiber) of $f$ over this point is the pullback
  $ X := {*} \times_Z Y$
  $$
    \raisebox{20pt}{
    \xymatrix{
        X \ar[r] \ar[d]& {*} \ar[d]
        \\
        Y \ar[r]^f & Z\,.
    }
	}
  $$
\end{definition}
\begin{observation}
  Let $f\colon Y\to Z$ in $\mathbf{H}$
  be as above.  Suppose that $Y$ is pointed and $f$ is a morphism of pointed objects.
  Then the $\infty$-fiber of an $\infty$-fiber is the loop object of the base.
\end{observation}
This means that we have a diagram
  $$
    \xymatrix{
        \Omega Z  \ar[d] \ar[r] & X \ar[r] \ar[d]& {*} \ar[d]
        \\
        {*} \ar[r] & Y \ar[r]^f & Z
    }
  $$
where the outer rectangle is an $\infty$-pullback if the left square is an
$\infty$-pullback. This follows from the pasting law, Proposition~\ref{PastingLawForPullbacks}.

\subsection{Epimorphisms and monomorphisms}
\label{StrucEpi}

In an $\infty$-topos there is an infinite tower of notions of epimorphisms and monomorphisms:
the $n$-connected and $n$-truncated morphisms for all $-2 \leq n \leq \infty$
\cite{Rezk, Lurie}.
\begin{definition}
 For $n \in \mathbb{N}$ an $\infty$-groupoid is called \emph{$n$-truncated} (or: an \emph{$n$-type})
 if all its homotopy groups in degree greater than $n$ are trivial.
 It is called \emph{$(-1)$-truncated} if it is either empty or contractible and
 \emph{$(-2)$-truncated} if it is non-empty and contractible. An object in an
 arbitrary $\infty$-category is $n$-truncated for $-2 \leq n < \infty$ if all
 hom-$\infty$-groupoids into it are $n$-truncated. A morphism of an $\infty$-category is called
 $n$-truncated if it is so as an object in the slice over its codomain
 (which means internally that its homotopy fibers are $n$-truncated).
 A $(-1)$-truncated morphism is also called a \emph{monomorphism}.
 The full embedding of the $n$-truncated objects
 of an $\infty$-topos is reflective, and the reflector $\tau_{\leq n}$ is called
 the \emph{$n$-truncation} operation.
\end{definition}
This is the topic of section 5.5.6 in \cite{Lurie}.
\begin{remark}
In a general $\infty$-topos every object has (groups of) \emph{homotopy sheaves} generalizing the
homotopy groups for bare $\infty$-groupoids.
If one knows that an object $X$ in an $\infty$-topos is truncated at all (for some possibly large truncation degree)
then it is $n$-truncated if all its homotopy sheaves $\pi_k(X)$ vanish in degree $k > n$.
\end{remark}
This is the content of Proposition 6.5.1.7 in \cite{Lurie}.
\begin{definition}
Let $\mathbf{H}$ be an $\infty$-topos.
For $X \to Y$ any morphism in $\mathbf{H}$, there is a
  simplicial object $\check{C}(X \to Y)$ in $\mathbf{H}$ (the {\em \v{C}ech
  nerve} of $f\colon X\to Y$) which in degree $n$ is the
  $(n+1)$-fold $\infty$-fiber product of $X$ over $Y$ with itself
  $$
    \check{C}(X \to Y) : [n] \mapsto X^{\times^{n+1}_Y}
  $$
  A morphism $f : X \to Y$ in $\mathbf{H}$ is an \emph{effective epimorphism}
  if it is the colimit of its own {\v C}ech nerve
  (under the natural map from the Cech nerve to $Y$):
  $$
    f : X \to \varinjlim \check{C}(X\to Y)
	\,.
  $$
  Write $\mathrm{Epi}(\mathbf{H}) \subset \mathbf{H}^I$ for the collection of
  effective epimorphisms.
  \label{EffectiveEpi}
\end{definition}
\begin{definition}
  For $n \in \mathbb{N}$ a morphism in an $\infty$-topos is called \emph{$n$-connected}
  if it is an effective epimorphism and all its homotopy sheaves are trivial in degree
  greater than $n$ when it is regarded as an object in the slice $\infty$-topos over its codomain.
  Any effective epimorphism is called $(-1)$-connected. An object $X$ is called $n$-connected
  if the canonical morphism $X \to \ast$ is $n$-connected.
\end{definition}
This is the topic of section 6.5.1 in \cite{Lurie}.
\begin{proposition}
  A morphism $f : X \to Y$ in the $\infty$-topos $\mathbf{H}$
  is an effective epimorphism precisely if its 0-truncation
  $\tau_0 f : \tau_0 X \to \tau_0 Y$ is an epimorphism (necessarily effective)
  in the 1-topos $\tau_{\leq 0} \mathbf{H}$.
  \label{EffectiveEpiIsEpiOn0Truncation}
\end{proposition}
This is Proposition 7.2.1.14 in \cite{Lurie}.
\begin{proposition}
  The classes $( \mathrm{Epi}(\mathbf{H}), \mathrm{Mono}(\mathbf{H}) )$
  constitute an orthogonal factorization system.
  \label{EpiMonoFactorizationSystem}
\end{proposition}
This is Proposition 8.5 in \cite{Rezk} and Example 5.2.8.16 in \cite{Lurie}.
\begin{definition}
  For $f : X \to Y$ a morphism in $\mathbf{H}$, we write its
  epi/mono factorization given by Proposition \ref{EpiMonoFactorizationSystem}
  as
  $$
    f :
    \xymatrix{
	  X \ar@{->>}[r] & \mathrm{im}(f)\  \ar@{^{(}->}[r] & Y
	}
  $$
  and we call $\xymatrix{\mathrm{im}(f)\ \ar@{^{(}->}[r] & Y}$ the \emph{$\infty$-image}
  of $f$.
  \label{image}
\end{definition}

\subsection{Groupoids and Groups}
\label{StrucInftyGroupoids}
\label{StrucInftyGroups}

In any $\infty$-topos $\mathbf{H}$ we may consider groupoids \emph{internal}
to $\mathbf{H}$, in the sense of internal category theory
(as exposed for instance in the introduction of \cite{Lurie2}).

Such a \emph{groupoid object}
$\mathcal{G}$ in $\mathbf{H}$ is an $\mathbf{H}$-object $\mathcal{G}_0$ ``of $\mathcal{G}$-objects''
together with an $\mathbf{H}$-object $\mathcal{G}_1$ ``of $\mathcal{G}$-morphisms''
equipped with source and target assigning morphisms $s,t : \mathcal{G}_1 \to \mathcal{G}_0$,
an identity-assigning morphism $i : \mathcal{G}_0 \to \mathcal{G}_1$ and a composition
morphism $\mathcal{G}_1 \times_{\mathcal{G}_0} \mathcal{G}_1 \to \mathcal{G}_1$
which together satisfy all the axioms of a groupoid (unitality, associativity, existence of
inverses) up to coherent homotopy in $\mathbf{H}$. One way to formalize what it
means for these axioms to hold up to coherent homotopy is as follows.

One notes that
ordinary groupoids, i.e.\  groupoid objects internal to $\mathrm{Set}$, are
characterized by the fact that their nerves are simplicial sets
$\mathcal{G}_\bullet : \Delta^{\mathrm{op}} \to \mathrm{Set}$
with the property that the groupoidal Segal maps
\[
\mathcal{G}_n\to \mathcal{G}_1\times_{\mathcal{G}_0}
\mathcal{G}_1\times_{\mathcal{G}_0} \cdots \times_{\mathcal{G}_0}
\mathcal{G}_1
\]
are isomorphisms for all $n\geq 2$.  This last condition
is stated precisely in Definition~\ref{GroupoidObject} below,
and clearly gives a characterization of groupoids that makes sense more generally, in
particular it makes sense internally to higher categories:
a groupoid object in $\mathbf{H}$ is an $\infty$-functor
$\mathcal{G} : \Delta^{\mathrm{op}} \to \mathbf{H}$ such that all groupoidal
Segal morphisms are equivalences in $\mathbf{H}$.
These $\infty$-functors $\mathcal{G}$ form the
objects of an $\infty$-category $\mathrm{Grpd}(\mathbf{H})$
of groupoid objects in $\mathbf{H}$.

Here a subtlety arises that is the source of a lot of interesting structure
in higher topos theory: the objects of $\mathbf{H}$ are themselves
 ``structured $\infty$-groupoids''. Indeed, there is a
full embedding $\mathrm{const} : \mathbf{H} \hookrightarrow \mathrm{Grpd}(\mathbf{H})$
that forms constant simplicial objects and thus regards every object $X \in \mathbf{H}$
as a groupoid object which, even though it has a trivial object of morphisms, already
has a structured $\infty$-groupoid of objects. This embedding is in fact
reflective, with the reflector given by forming the $\infty$-colimit
over a simplicial diagram, the ``geometric realization''
$$
  \xymatrix{
    \mathbf{H}
      \ar@{<-}@<+1.25ex>[rr]^-{\varinjlim}
	  \ar@{^{(}->}@<-1.25ex>[rr]_-{\mathrm{const}}^-{\perp}
	  &&
	  \mathrm{Grpd}(\mathbf{H})
  }
  \,.
$$
For $\mathcal{G}$ a groupoid object in $\mathbf{H}$, the object
$\varinjlim \mathcal{G}_\bullet$ in $\mathbf{H}$
may be thought of as the
$\infty$-groupoid obtained by ``gluing together the object of objects of
$\mathcal{G}$ along the object of morphisms of $\mathcal{G}$''.
This idea that groupoid objects in an $\infty$-topos are
like structured $\infty$-groupoids together with gluing information
is formalized by the statement recalled as
Theorem~\ref{NaturalThirdGiraud} below, which says that groupoid objects in
$\mathbf{H}$ are equivalent to the \emph{effective epimorphisms}
$\xymatrix{Y \ar@{->>}[r] & X}$ in $\mathbf{H}$, the intrinsic notion of
\emph{cover} (of $X$ by $Y$) in $\mathbf{H}$. The effective epimorphism/cover
corresponding to a groupoid object $\mathcal{G}$ is the colimiting cocone
$\xymatrix{\mathcal{G}_0 \ar@{->>}[r] & \varinjlim \mathcal{G}_\bullet}$.

\medskip

After this preliminary discussion we state the
following definition of groupoid object in
an $\infty$-topos (this definition appears in \cite{Lurie}
as Definition 6.1.2.7, using Proposition 6.1.2.6).

\begin{definition}[\cite{Lurie}, Definition 6.1.2.7]
  \label{GroupoidObject}
  A \emph{groupoid object} in an $\infty$-topos $\mathbf{H}$ is
  a simplicial object
  $$
    \mathcal{G} : \Delta^{\mathrm{op}} \to \mathbf{H}
  $$
  all of whose groupoidal Segal maps are equivalences:
  in other words, for every
  $n \in \mathbb{N}$
  and every partition $[k] \cup [k']  =  [n]$ into two subsets
  such that $[k] \cap [k'] = \{*\}$, the canonical diagram
  \[
    \xymatrix{
      \mathcal{G}_n \ar[r] \ar[d] & \mathcal{G}_k \ar[d]
      \\
      \mathcal{G}_{k'} \ar[r] & \mathcal{G}_{0}
    }
  \]
  is an $\infty$-pullback diagram.  We write
  \[
    \mathrm{Grpd}(\mathbf{H}) \subset \mathrm{Func}(\Delta^{\mathrm{op}}, \mathbf{H})
  \]
  for the full subcategory of the $\infty$-category of simplicial
  objects in $\mathbf{H}$ on the groupoid objects.
\end{definition}
The following example is fundamental. In fact the third $\infty$-Giraud axiom
says that up to equivalence, all groupoid objects are of this form.
\begin{example}
  For $X \to Y$ any morphism in $\mathbf{H}$, the \v{C}ech
  nerve $\check{C}(X\to Y)$ of $X\to Y$ (Definition~\ref{EffectiveEpi}) is a
  groupoid object.  This example appears in \cite{Lurie} as Proposition 6.1.2.11.
\end{example}

The following statement refines
the third $\infty$-Giraud axiom, Definition \ref{GiraudRezkLurieAxioms}.
\begin{theorem}
  \label{NaturalThirdGiraud}
   There is a natural equivalence of $\infty$-categories
  $$
    \mathrm{Grpd}(\mathbf{H})
     \simeq
    (\mathbf{H}^{\Delta[1]})_{\mathrm{eff}}
    \,,
  $$
  where $(\mathbf{H}^{\Delta[1]})_{\mathrm{eff}}$
  is the full sub-$\infty$-category of the
  arrow category $\mathbf{H}^{\Delta[1]}$
  of $\mathbf{H}$ on the effective epimorphisms, Definition \ref{EffectiveEpi}.
\end{theorem}
This appears below Corollary 6.2.3.5 in \cite{Lurie}.

\medskip

In addition, every $\infty$-topos $\mathbf{H}$
comes with a notion of \emph{$\infty$-group objects} that generalize both the
ordinary notion of group objects in a topos as well as that of
grouplike $A_\infty$-spaces in $\mathrm{Grpd}_{\infty}$.

\medskip

There is an evident definition of what an $\infty$-group object in $\mathbf{H}$
should be, and then there is a theorem saying that this is equivalent to a certain kind of
simplicial object in $\mathbf{H}$.
This theorem is part of what, we find, makes the theory of groups, group actions and principal
bundles in an $\infty$-topos be so well behaved, and we will mostly work with this simplicial
incarnation of group objects. But the evident definition that the reasoning starts with is
of course this:
a group object is an object which is equipped with an associative and unital
product operation such that for each element there is an inverse. Now in the
homotopy-theoretic context of $\infty$-topos theory an associative unital structure
means an associative unital structure \emph{up to coherent homotopy} and the technical
term for this is \emph{$A_\infty$-structure}, famous from the theory of
loop spaces, see \cite{LurieAlgebra} for a comprehensive
modern account. Moreover, statements about elements here are supposed to be statements
about connected components, and hence we ask for such $A_\infty$ structures
such that on connected components the product operation is invertible (such
$A_\infty$-structures are traditionally also called ``groupal'' or ``grouplike'').

Therefore the manifest definition of $\infty$-group objects in $\mathbf{H}$ is the following
(this appears as Definition 5.1.3.2 together with
Remark 5.1.3.3 in \cite{LurieAlgebra}).
\begin{definition}
\label{inftygroupinootopos}
  An \emph{$\infty$-group} in $\mathbf{H}$
  is an $A_\infty$-algebra $G$ in $\mathbf{H}$ such that
  the sheaf of connected components
  $\pi_0(G)$ is a group object in $\tau_{\leq 0} \mathbf{H}$.
  Write
  $\mathrm{Grp}(\mathbf{H})$ for the $\infty$-category
  of $\infty$-groups in $\mathbf{H}$.
\end{definition}
As in classical algebraic topology, the fundamental examples of such $\infty$-groups
arise from forming loops, and there is a central de-looping theorem saying that, up to equivalence,
in fact all $\infty$-groups arise this way:
\begin{definition}
   Write
  \begin{itemize}
    \item $\mathbf{H}^{*/}$ for the
  $\infty$-category of pointed objects in $\mathbf{H}$;
    \item  $\mathbf{H}_{\geq 1}$
    for the full sub-$\infty$-category of $\mathbf{H}$ on the
    connected objects;
    \item
      $\mathbf{H}^{*/}_{\geq 1}$ for the full sub-$\infty$-category
      of the pointed objects on the connected objects.
   \end{itemize}
\end{definition}
\begin{definition}
  \label{loop space object}
  Write
  $$
    \Omega : \mathbf{H}^{*/} \to \mathbf{H}
  $$
  for the $\infty$-functor that sends a pointed object $* \to X$
  to its \emph{loop space object}, i.e.\  the $\infty$-pullback
  $$
   \raisebox{20pt}{
   \xymatrix{
      \Omega X \ar[r]\ar[d] & {*} \ar[d]
      \\
      {*} \ar[r] & X\,.
   }
   }
  $$
\end{definition}
\begin{theorem}[Lurie]
  \label{DeloopingTheorem}
  \label{delooping}
  Every loop space object canonically has the structure of an
  $\infty$-group, and this construction extends to an
  $\infty$-functor
  $$
    \Omega : \mathbf{H}^{*/} \to \mathrm{Grp}(\mathbf{H})
    \,.
  $$
  This $\infty$-functor constitutes part of an equivalence of $\infty$-categories
  $$
    (\Omega \dashv \mathbf{B})
    :
    \xymatrix{
      \mathrm{Grp}(\mathbf{H})
       \ar@{<-}@<+5pt>[r]^<<<<<<{\Omega}
       \ar@<-5pt>[r]_<<<<<<{\mathbf{B}}^<<<<<<\simeq
      &
      \mathbf{H}^{*/}_{\geq 1}\, .
    }
  $$
\end{theorem}
This is Lemma 7.2.2.1 in \cite{Lurie}.
(See also Theorem 5.1.3.6 of \cite{LurieAlgebra}
where this is the equivalence denoted $\phi_0$ in the proof.) For
$\mathbf{H} = \mathrm{Grpd}_{\infty}$ this
reduces to various classical
theorems in homotopy theory, for instance the
construction of classifying spaces (Kan and Milnor) and de-looping theorems (May and Segal).
\begin{definition}
We call the inverse
$\mathbf{B} : \mathrm{Grp}(\mathbf{H}) \to \mathbf{H}^{*/}_{\geq 1}$
in Theorem~\ref{DeloopingTheorem} above the
\emph{delooping} functor of $\mathbf{H}$. By convenient abuse
of notation we write $\mathbf{B}$ also for the composite
$\mathbf{B} : \mathrm{Grpd}(\mathbf{H}) \to \mathbf{H}^{*/}_{\geq 1}
\to \mathbf{H}$ with the functor that forgets the basepoint and the
connectivity.
\end{definition}
\begin{remark}
  Even if the connected objects involved admit an essentially
  unique point,
  the homotopy type of the full hom-$\infty$-groupoid
  $\mathbf{H}^{*/}(\mathbf{B}G, \mathbf{B}H)$
  of pointed objects in general differs
  from the hom $\infty$-groupoid $\mathbf{H}(\mathbf{B}G, \mathbf{B}H)$
  of the underlying unpointed objects.
  For instance let $\mathbf{H} := \mathrm{Grpd}_{\infty}$ and let $G$ be
  an ordinary group, regarded as a group object in $\mathrm{Grpd}_{\infty}$.
  Then the invertible elements in $\mathbf{H}^{*/}(\mathbf{B}G, \mathbf{B}G)$
  give the ordinary automorphism group $\mathrm{Aut}(G)$ of $G$, but
  the invertible elements in $\mathbf{H}(\mathbf{B}G, \mathbf{B}G)$ is the
  automorphism 2-group of $G$, we discuss this
  further around Example~\ref{automorphism2GroupAbstractly} below.
\end{remark}
Now observe that for $X$ a pointed connected object, then the point inclusion $\ast \to X$
is an effective epimorphism and the loop space object $\Omega X$ in def. \ref{loop space object}
is the first stage of the corresponding {\v C}ech nerve, as in the discussion of
groupoid objects above in \ref{StrucInftyGroupoids}.
This suggests that, moreover, group objects in $\mathbf{H}$ should be equivalent
to those groupoid objects whose degree-0 piece is equivalent to the point. This is
indeed the case, and this is central to the development of our discussion:
\begin{proposition}[Lurie]
  \label{InfinityGroupObjectsAsGroupoidObjects}
  $\infty$-groups $G$ in $\mathbf{H}$ are equivalently
  those groupoid objects $\mathcal{G}$ in $\mathbf{H}$ (Definition~\ref{GroupoidObject})
  for which $\mathcal{G}_0 \simeq *$.
\end{proposition}
This is the statement of the compound equivalence
$\phi_3\phi_2\phi_1$ in the proof of Theorem 5.1.3.6 in
\cite{LurieAlgebra}.
\begin{remark}
  \label{PointIntoBGIsEffectiveEpimorphism}
  \label{Cech nerve of * -> BG}
  This means that for $G$ an $\infty$-group object,
  the {\v C}ech nerve extension of its delooping fiber
  sequence $G \to * \to \mathbf{B}G$ is the simplicial
  object
  $$
    \xymatrix{
       \cdots
       \ar@<+6pt>[r] \ar@<+2pt>[r] \ar@<-2pt>[r] \ar@<-6pt>[r]
       &
       G \times G
       \ar@<+4pt>[r]
       \ar[r]
        \ar@<-4pt>[r]
         &
         G
         \ar@<+2pt>[r] \ar@<-2pt>[r]
       &
       {*}
       \ar@{->>}[r]
       &
       \mathbf{B}G
    }
  $$
  that exhibits $G$ as a groupoid object over $*$.
  In particular it means that for $G$ an $\infty$-group, the
%
%
%
%
given
morphism $* \to \mathbf{B}G$
  is an effective epimorphism.
  \end{remark}

\subsection{Cohomology}
\label{StrucCohomology}
\label{section.Cohomology}

There is an intrinsic notion of \emph{cohomology}
in every $\infty$-topos $\mathbf{H}$: it is simply given by the
connected components of mapping spaces. Of course such mapping spaces
exist in every $\infty$-category, but we need some extra conditions on
$\mathbf{H}$ in order for them to behave like cohomology sets. For
instance, if $\mathbf{H}$ has pullbacks then there is a notion of
long exact sequences in cohomology.
Our main theorem (Theorem~\ref{PrincipalInfinityBundleClassification} below) will show that the second
and third $\infty$-Giraud axioms imply that this intrinsic
notion of cohomology has the property that it \emph{classifies}
certain geometric structures in the $\infty$-topos.

\begin{definition}
 \label{cohomology}
For $X,A \in \mathbf{H}$ two objects, we say that
$$
  H^0(X,A) := \pi_0 \mathbf{H}(X,A)
$$
is the \emph{cohomology set}\index{cohomology!general abstract} of $X$ with coefficients in $A$.
In particular if $G$ is an  $\infty$-group we write
$$
  H^1(X,G) := H^0(X,\mathbf{B}G) = \pi_0 \mathbf{H}(X, \mathbf{B}G)
$$
for cohomology with coefficients in the delooping
$\mathbf{B}G$ of $G$.
Generally, if  $K \in \mathbf{H}$ has a
specified
$n$-fold delooping
$\mathbf{B}^nK$ for some non-negative integer $n$, we write
$$
  H^n(X,K) := H^0(X,\mathbf{B}^n K) = \pi_0 \mathbf{H}(X, \mathbf{B}^n K)
  \,.
$$
\end{definition}
In the context of cohomology on $X$ with coefficients in $A$ we say that
\begin{itemize}
\item
 the hom-space $\mathbf{H}(X,A)$ is the \emph{cocycle $\infty$-groupoid}\index{
!cocycle};
\item
  an object $g : X \to A$ in $\mathbf{H}(X,A)$ is a \emph{cocycle};
\item
  a morphism: $g \Rightarrow h$ in $\mathbf{H}(X,A)$ is a \emph{coboundary} between cocycles.
\item
  a morphism $c : A \to B$ in $\mathbf{H}$
  represents the \emph{universal characteristic class}\index{characteristic class!general abstract}
  (cohomology operation)
  $$
    [c] : H^0(-,A) \to H^0(-,B)
    \,.
  $$
\end{itemize}
If $X\simeq Y/\!/G$ is a homotopy quotient, then the cohomology of $X$ is
equivariant cohomology of $Y$. Similarly, for general $X$ this notion of cohomology
incorporates various local notions of equivariance (for instance $X$ might be an orbifold
which is only locally equivalent to a global quotient).
\begin{remark}
 \label{CohomologyOverX}
Of special interest is the cohomology defined by a slice
$\infty$-topos
$$
  \mathcal{X} := \mathbf{H}_{/X}
$$
over some $X \in \mathbf{H}$.
Such a slice is canonically equipped with the
{\'e}tale geometric morphism (\cite{Lurie}, Remark 6.3.5.10)
$$
  ((p_X)_! \dashv (p_X)^* \dashv (p_X)_*)
  :
  \xymatrix{
     \mathbf{H}_{/X}
	 \ar@<1.3ex>[rr]^-{(p_X)_!}
	 \ar@{<-}[rr]|-{(p_X)^*}
	 \ar@<-1.3ex>[rr]_-{(p_X)_*}
	 &&
	 \mathbf{H}
  }
  \,,
$$
where $p_X : X\to \ast$ is the canonical morphism,
$(p_X)_!$ simply forgets the morphism to $X$ and
where $(p_X)^* = X \times (-)$ forms the product with $X$.
Accordingly we have
$(p_X)^* (*_{\mathbf{H}}) \simeq *_{\mathcal{X}}$
and $(p_X)_! (*_{\mathcal{X}}) = X \in \mathbf{H}$,
saying that the terminal object $\ast_{\mathcal{X}}$ in $\mathcal{X}$, which
is the pullback of the terminal object $\ast_{\mathbf{H}}$ of $\mathbf{H}$ to $\mathcal{X}$,
is identified with $X$ itself.
Therefore
cohomology over $X$ with coefficients of the form $(p_X)^* A$ is
equivalently the cohomology in $\mathbf{H}$ of $X$ with coefficients in $A$:
$$
  \mathcal{X}(X, (p_X)^* A) \simeq \mathbf{H}(X,A)
  \,.
$$
But for a general coefficient object $A \in \mathcal{X}$ the
$A$-cohomology over $X$ in $\mathcal{X}$ is a
\emph{twisted} cohomology of $X$ in $\mathbf{H}$.
This we discuss below in Section~\ref{TwistedCohomology}.
\end{remark}
Typically one thinks of a morphism $A \to B$ in
$\mathbf{H}$ as presenting a \emph{characteristic class} of $A$ if
$B$ is ``simpler'' than $A$, notably if $B$ is an Eilenberg-MacLane object $B = \mathbf{B}^n K$ for
$K$ a 0-truncated abelian group in $\mathbf{H}$. In this case the characteristic class may
be regarded as being in the degree-$n$ $K$-cohomology of $A$
$$
  [c] \in H^n(A,K)
  \,.
$$
\begin{definition}
 \label{fiber sequence}
 \label{LongFiberSequence}
For every morphism $c : \mathbf{B}G \to \mathbf{B}H \in \mathbf{H}$ define the
\emph{long fiber sequence to the left}
$$
   \cdots
   \to
   \Omega G
  \to
  \Omega H
   \to
   \Omega F
    \to
   G
     \to
   H
     \to
    F
     \to
    \mathbf{B}G
      \stackrel{c}{\to}
   \mathbf{B}H
$$
by the consecutive pasting diagrams of $\infty$-pullbacks
$$
  \xymatrix{
      & \ar@{..}[d] & \ar@{..}[d]
      \\
     \ar@{..}[r] &  \Omega F \ar[d]\ar[r]& G \ar[r] \ar[d] & {*} \ar[d]
      \\
      & {*} \ar[r] & H \ar[r] \ar[d] & F \ar[r] \ar[d] & {*} \ar[d]
      \\
      & & {*} \ar[r] & \mathbf{B}G \ar[r]^c & \mathbf{B}H
  }
$$
\end{definition}
We have the following basic fact.
\begin{theorem}
\begin{enumerate}
\item In the long fiber sequence to the left of
$c : \mathbf{B}G \to \mathbf{B}H$ after $n$ iterations
all terms are equivalent to the point if $H$ and $G$ are $n$-truncated.
\item For every object $X \in \mathbf{H}$ we have a long exact sequence of pointed cohomology sets
  $$
    \cdots \to H^0(X,G) \to H^0(X,H) \to H^1(X,F) \to H^1(X,G) \to H^1(X,H)
   \,.
  $$
\end{enumerate}
 \label{LongExactSequenceInCohomology}
\end{theorem}
%
\begin{remark}
  For the special case that $G$ is a 1-truncated $\infty$-group (or
  \emph{2-group}) classified by a 3-cocycle $\mathbf{c}$,
  Theorem \ref{LongExactSequenceInCohomology} is
  a classical result due to \cite{BreenBitorseurs}. The first and
  only nontrivial
  stage of the internal Postnikov tower
  $$
    \xymatrix{
	  \mathbf{B}^2 A \ar[r] & \mathbf{B}G \ar[d]
	  \\
	  & \mathbf{B} H
      \ar[r]^-{\mathbf{c}}
      &
      \mathbf{B}^3 A
	}
  $$
  of the delooped 2-group (with $H := \tau_0 G\in \tau_{\leq 0} \mathrm{Grp}(\mathbf{H})$
  an ordinary group object and $A := \pi_1 G \in \tau_{\leq 0} \mathrm{Grp}(\mathbf{H})$
  an ordinary abelian group object) yields the long exact sequence of pointed cohomology
  sets
  $$
    0 \to H^1(-,A) \to H^0(-,G) \to H^0(-,H) \to H^2(-,A) \to H^1(-,G) \to
	H^1(-,H) \to H^3(-,A)
  $$
  (see also \cite{NikolausWaldorf2}.)
  Notably, the last morphism gives the obstructions against lifting traditional
  nonabelian cohomology $H^1(-,H)$ to nonabelian cohomology $H^1(-,G)$ with values
  in the 2-group. This we discuss further in Section \ref{ExtensionsOfCohesiveInfinityGroups}.
  \label{ReferencesOnLongSequences}
\end{remark}

Generally, to every cocycle $g : X \to \mathbf{B}G$ is
canonically associated its $\infty$-fiber
$P \to X$ in $\mathbf{H}$, the $\infty$-pullback
$$
  \raisebox{20pt}{
  \xymatrix{
    P \ar[r] \ar[d]& {*} \ar[d]
    \\
    X \ar[r] ^g & \mathbf{B}G
    \,.
  }
  }
$$
We now discuss how each such $P$ canonically has the structure of a
\emph{$G$-principal $\infty$-bundle} and that $\mathbf{B}G$ is the
\emph{fine moduli object} (the \emph{moduli $\infty$-stack}) for
$G$-principal $\infty$-bundles.\footnote{The concept of (fine) moduli stacks
is historically most commonly associated with algebraic geometry, but the problem which they solve, namely the
classification of structures including their (auto-)equivalences, is universal.
Specifically, if $\mathbf{H}$ is the $\infty$-topos over a site of schemes then it
contains the moduli stacks as they appear in algebraic geometry. }

\section{Principal bundles}
\label{Principal infinity-bundles general abstract}
\label{PrincipalInfBundle}

We define here $G$-principal $\infty$-bundles in any $\infty$-topos
$\mathbf{H}$, discuss their basic properties and show that they are classified
by the intrinsic $G$-cohomology in $\mathbf{H}$, as discussed in
Definition~\ref{cohomology}.

\subsection{Introduction and survey}
\label{PrincBund-Intro}

Let $G$ be a topological group, or Lie group or
some similar such object. The traditional
definition of \emph{$G$-principal bundle} is the following:
there is a map
$$
  P \to X := P/G
$$
which is the quotient projection
induced by a \emph{free} action
$$
  \rho : P \times G \to P
$$
of $G$ on a space (or manifold, depending on context) $P$,
such that there is a cover $U \to X$ over which the quotient projection is isomorphic
to the trivial one $U \times G \to U$.

In higher geometry, if $G$ is a topological or smooth
$\infty$-group, the quotient projection must be
replaced by the $\infty$-quotient (homotopy quotient)
projection
\[
P\to X := P/\!/ G
\]
for the action of $G$ on a topological or smooth $\infty$-groupoid
(or $\infty$-stack) $P$.  It is a remarkable fact that this
single condition on the map $P\to X$
already implies that $G$ acts freely on $P$ and that $P\to X$
is locally trivial, when the latter notions are understood in the
context of higher geometry.  We will therefore define
a $G$-principal $\infty$-bundle to be such a map $P\to X$.


As motivation for this, notice that if a Lie group $G$ acts properly,
but not freely, then the quotient $P \to X := P/G$ differs from the homotopy quotient,
which instead is locally the quotient stack by the non-free part of the group
action (an orbifold, if the stabilizers are finite).

Conversely this means that in the context of higher geometry a non-free action
may also be principal: with respect not to a base space, but with respect to a base groupoid/stack.
In the example just discussed, we have that the projection $P \to X/\!/ G_{\mathrm{stab}}$
exhibits $P$ as a $G$-principal bundle over the action groupoid
$P /\!/ G \simeq X/\!/ G_{\mathrm{stab}}$. For instance if $P = V$ is a
vector space equipped with a $G$-representation, then $V \to V/\!/ G$ is a
$G$-principal bundle over a groupoid/stack.
In other words, the traditional requirement of freeness in a principal action is not so much
a characterization of principality as such, as rather a condition that ensures that the
base of a principal action is a 0-truncated object in higher geometry.

Beyond this specific class of 0-truncated examples, this means that we have the following
noteworthy general statement: in higher geometry \emph{every} $\infty$-action
is principal with respect to
\emph{some} base, namely with respect to its $\infty$-quotient.

More is true: in the context of an $\infty$-topos
every $\infty$-quotient projection of an $\infty$-group action
is locally trivial, with respect to
the canonical intrinsic notion of cover, hence of locality. Therefore
also the condition of local triviality in the classical definition of principality
becomes automatic.  In fact, from the $\infty$-Giraud axioms, we
see that the projection map $P \to P /\!/ G$ is always a cover
(an \emph{effective epimorphism}) and so, since every $G$-principal $\infty$-bundle
trivializes over itself, it exhibits a local trivialization of itself;
even without explicitly requiring it to be locally trivial.

As before, this means that the local triviality clause appearing in the
traditional definition of principal bundles is not so much a characteristic of
principality as such, as rather a condition that ensures that a given quotient
taken in a category of geometric spaces coincides with the ``refined'' quotient
obtained when regarding the situation in the ambient $\infty$-topos.

Another direct consequence of the $\infty$-Giraud axioms
is the equivalence of the definition of principal bundles as quotient maps,
as we  have discussed so far, with the other main definition of principality: the condition
that the ``shear map'' $ (\mathrm{id}, \rho) :  P \times G \to P \times_X P$ is an equivalence.
It is immediate to verify in traditional 1-categorical contexts that this is
equivalent to the action being properly free and exhibiting $X$ as its quotient.
Simple as this is, one may observe, in view of the above discussion,
that the shear map being an equivalence is much more fundamental even: notice
that $P \times G$ is the first stage of the \emph{action groupoid object}
$(P/\!/G)_\bullet$, and that $P \times_X P$ is the first stage of the \emph{{\v C}ech nerve groupoid object}
$\check{C}(P \to X)$ of the corresponding quotient map. Accordingly, the shear map equivalence
is the first stage in the equivalence of groupoid objects in the $\infty$-topos
$$
  (P /\!/G)_\bullet \simeq \check{C}(P \to X)
  \,.
$$
This equivalence is just the explicit statement of the fact mentioned before: the groupoid object
$(P/\!/G)_\bullet$ is effective -- as is any groupoid object in an $\infty$-topos -- and, equivalently,
its principal $\infty$-bundle map $P \to X$ is an effective epimorphism.

Fairly directly from this fact, finally, springs the classification theorem of
principal $\infty$-bundles. For we have a canonical morphism of groupoid objects
$(P /\!/G)_\bullet \to (* / \!/G)_\bullet$ induced by the terminal map $P \to *$. By the $\infty$-Giraud
theorem the $\infty$-colimit over this sequence of morphisms of groupoid objects
is a $G$-cocycle on $X$ (Definition~\ref{cohomology}) canonically induced by $P$:
$$
  \varinjlim  \left(\check{C}(P \to X)_\bullet \simeq (P /\!/G)_\bullet \to (* /\!/G)_\bullet \right)
    =
  (X \to \mathbf{B}G)
  \;\;\;
  \in \mathbf{H}(X, \mathbf{B}G)
  \,.
$$
Conversely, from any such $G$-cocycle one obtains a $G$-principal
$\infty$-bundle simply by forming its $\infty$-fiber: the $\infty$-pullback of
the point inclusion ${*} \to \mathbf{B}G$.  We show in \cite{NSSb} that in presentations
of the $\infty$-topos theory by 1-categorical tools, the computation of this homotopy
fiber is \emph{presented} by the ordinary pullback of a big resolution of the point,
which turns out to be nothing but the universal $G$-principal bundle.
This appearance of the universal $\infty$-bundle as just
a resolution of the point inclusion may be understood in light of the above discussion
as follows.
The classical characterization of the
universal $G$-principal bundle $\mathbf{E}G$ is as a space that is homotopy equivalent
to the point and equipped with a \emph{free} $G$-action. But by the above, freeness of the
action is an artefact of 0-truncation and not a characteristic of principality in higher
geometry. Accordingly, in higher geometry the universal $G$-principal $\infty$-bundle
for any $\infty$-group $G$ may be taken to \emph{be} the point, equipped with the
trivial (maximally non-free) $G$-action. As such, it is a bundle not over the
classifying \emph{space} $B G$ of $G$, but over the full moduli $\infty$-stack $\mathbf{B}G$.

The following table summarizes the relation between
$\infty$-bundle theory and the $\infty$-Giraud axioms as indicated above, and as
proven in the following section.

\medskip
\begin{center}
\begin{tabular}{c|c}
 {\bf $\infty$-Giraud axioms} & {\bf principal $\infty$-bundle theory}
 \\
 \hline
 \hline
 quotients are effective &
   \begin{tabular}{c} \\ every $\infty$-quotient $P \to X := P/\!/ G$ \\is principal \\ \ \end{tabular}
 \\
 \hline
  colimits are preserved by pullback &
  \begin{tabular}{c}\\ $G$-principal $\infty$-bundles \\ are classified by $\mathbf{H}(X,\mathbf{B}G)$\\
   \end{tabular}
\end{tabular}
\end{center}

\subsection{Definition and classification}
\label{DefintionAndClassification}

\begin{definition}
  \label{ActionInPrincipal}
  For $G \in \mathrm{Grp}(\mathbf{H})$ a group object,
  we say a \emph{$G$-action} on an object $P \in \mathbf{H}$
  is a groupoid object $(P/\!/G)_\bullet$ (Definition~\ref{GroupoidObject}) of the form
  $$
    \xymatrix{
       \cdots
       \ar@<+6pt>[r] \ar@<+2pt>[r] \ar@<-2pt>[r] \ar@<-6pt>[r]
       &
       P \times G \times G
       \ar@<+4pt>[r]
       \ar[r]
        \ar@<-4pt>[r]
         &
         P \times G
         \ar@<+2pt>[r]^<<<<<{\rho := d_0 } \ar@<-2pt>[r]_{d_1}
       & P
    }
  $$
  such that $d_1 : P \times G \to P$ is the projection, and such that
  the degreewise projections
  $P \times G^n \to G^n $ constitute a morphism of groupoid
  objects
  $$
    \xymatrix{
       \cdots
       \ar@<+6pt>[r] \ar@<+2pt>[r] \ar@<-2pt>[r] \ar@<-6pt>[r]
       &
       P \times G \times G
       \ar[d]
       \ar@<+4pt>[r]
       \ar[r]
        \ar@<-4pt>[r]
         &
         P \times G
         \ar[d]
         \ar@<+2pt>[r] \ar@<-2pt>[r]
       & P
       \ar[d]
       \\
       \cdots
       \ar@<+6pt>[r] \ar@<+2pt>[r] \ar@<-2pt>[r] \ar@<-6pt>[r]
       &
       G \times G
       \ar@<+4pt>[r]
       \ar[r]
        \ar@<-4pt>[r]
         &
         G
         \ar@<+2pt>[r] \ar@<-2pt>[r]
       & {*}
    }
  $$
where the lower simplicial object exhibits $G$ as a groupoid
object over $\ast$ (see Remark~\ref{Cech nerve of * -> BG}).

  With convenient abuse of notation we also write
  $$
    P/\!/G := \varinjlim (P \times G^{\times^\bullet})\;\;
    \in \mathbf{H}
  $$
  for the corresponding $\infty$-colimit object, the \emph{$\infty$-quotient} of this
  action.

  Write
  $$
    G \mathrm{Action}(\mathbf{H}) \hookrightarrow \mathrm{Grpd}(\mathbf{H})_{/({*}/\!/G)}
  $$
  for the full sub-$\infty$-category of groupoid objects over $*/\!/G$
  having as objects those that are $G$-actions.
\end{definition}
\begin{remark}
  \label{ActionMapEncodedInGBundle}
  The remaining face map $d_0$ in Definition \ref{ActionInPrincipal}
  $$
    \rho := d_0 : P \times G \to P
  $$
  is the action itself and the condition that it fits into such a simplicial diagram
  encodes the action property up to coherent homotopy. For instance
  using effectivity of groupoid objects, Definition~\ref{GiraudRezkLurieAxioms},
  and the defining assumption on $(P/\!/G)_1$
  it follows that we have specified equivalences (where we abbreviate $X := P/\!/G$)
  $$
    (P/\!/G)_2
    \simeq
    (P \times_X P) \times_X P
       \simeq
    (P \times G) \times_X P
       \simeq
    P \times G \times G
    \,.
  $$
  From this it follows that the three maps from $P \times G \times G$ to $P \times G$ here are given by, respectively,
  multiplication of the two group factors, action of the first group factor on $P$ and projection on $P \times G$.
  The simplicial identities in $\mathbf{H}$ then give, in particular, a homotopy between, first, the composition
  of multiplying in the group and then acting and, second, the composition of acting with one factor
  and then with the other. At the next stage the simplicial identities encode that this homotopy in turn is compatible
  with the associativity-homotopy involved in acting with three group factors, and so on.
\end{remark}
\begin{remark}
  Using this notation in Proposition~\ref{InfinityGroupObjectsAsGroupoidObjects}
  we have
  $$
    \mathbf{B}G \simeq */\!/G
    \,.
  $$
\end{remark}
We list examples of $\infty$-actions below as Example \ref{ExamplesOfActions}.
This is most conveniently done after establishing the theory of
principal $\infty$-actions, to which we now turn.

\begin{definition}
  \label{principalbundle}
  Let $G \in \infty \mathrm{Grp}(\mathbf{H})$ be an
  $\infty$-group and let $X$ be an object of $\mathbf{H}$.
  A {\em $G$-principal $\infty$-bundle} over $X$
  (or \emph{$G$-torsor over $X$})
  is
  \begin{enumerate}
    \item a morphism $P \to X$ in $\mathbf{H}$;
	\item together with a $G$-action on $P$;
  \end{enumerate}
  such that $P \to X$ is the colimiting cocone exhibiting the quotient map
  $X \simeq P/\!/G$ (Definition \ref{ActionInPrincipal}).

  A \emph{morphism} of $G$-principal $\infty$-bundles over $X$ is a morphism of $G$-actions
  that fixes $X$; the $\infty$-category of $G$-principal $\infty$-bundles over $X$
  is the homotopy fiber of $\infty$-categories
  $$
    G \mathrm{Bund}(X) := G \mathrm{Action}(\mathbf{H}) \times_{\mathbf{H}} \{X\}
  $$
  over $X$ of the quotient map
  $$
    \xymatrix{
      G \mathrm{Action}(\mathbf{H})
	  \ar@{^{(}->}[r] & \mathrm{Grpd}(\mathbf{H})_{/(*/\!/G)}
	  \ar[r] &
	  \mathrm{Grpd}(\mathbf{H})
	  \ar[r]^-{\varinjlim}
	  &
	  \mathbf{H}
	}
	\,.
  $$
\end{definition}
\begin{remark}
  \label{GBundlesAreEffectiveEpimorphisms}
  By the third $\infty$-Giraud axiom (Definition~\ref{GiraudRezkLurieAxioms})
  this means in particular that a
  $G$-principal $\infty$-bundle $P \to X$ is an
  effective epimorphism in $\mathbf{H}$.
\end{remark}
\begin{remark}
  Even though $G \mathrm{Bund}(X)$ is by definition a priori an $\infty$-category,
  Proposition \ref{MorphismsOfInfinityBundlesAreEquivalences} below says
  that in fact it happens to be an $\infty$-groupoid: all
  its morphisms are invertible.
\end{remark}
\begin{proposition}
  \label{PrincipalityCondition}
  A $G$-principal $\infty$-bundle $P \to X$ satisfies the
  \emph{principality condition}: the canonical morphism
  $$
    (\rho, p_1)
	:
    \xymatrix{
	  P \times G
	  \ar[r]^{\simeq}
	  &
	  P \times_X P
	}
  $$
  is an equivalence, where $\rho$ is the $G$-action.
\end{proposition}
\proof
  By the third $\infty$-Giraud axiom (Definition~\ref{GiraudRezkLurieAxioms}) the groupoid object
  $P/\!/G$ is effective, which means that it is equivalent
  to the   {\v C}ech nerve of $P \to X$. In first degree this implies
  a canonical equivalence $P \times G \to P \times_X P$. Since
  the two face maps $d_0, d_1 : P \times_X P \to P$ in the
  {\v C}ech nerve are simply the projections out of the fiber product,
  it follows that the two components of this canonical equivalence
  are the two face maps $d_0, d_1 : P \times G \to P$ of $P/\!/G$.
  By definition, these are the projection onto the first factor
  and the action itself.
\endofproof
\begin{proposition}
  \index{principal $\infty$-bundle!construction from cocycle}
  \label{BundleStructureOnInfinityFiber}
  \label{PrincipalBundleAsHomotopyFiber}
  For $g : X \to \mathbf{B}G$ any morphism, its homotopy fiber
  $P \to X$ canonically carries the structure of a
  $G$-principal $\infty$-bundle over $X$.
\end{proposition}
\proof
  That $P \to X$ is the fiber of $g : X \to \mathbf{B}G$
  means that we have an $\infty$-pullback diagram
  $$
    \raisebox{20pt}{
    \xymatrix{
      P \ar[d]\ar[r] & {*}\ar[d]
      \\
      X \ar[r]^g & \mathbf{B}G.
    }
	}
  $$
  By the pasting law for $\infty$-pullbacks, Proposition~\ref{PastingLawForPullbacks},
  this induces a compound diagram
  $$
    \xymatrix{
       \cdots
       \ar@<+6pt>[r] \ar@<+2pt>[r] \ar@<-2pt>[r] \ar@<-6pt>[r]
       &
       P \times G \times G
       \ar[d]
       \ar@<+4pt>[r]
       \ar[r]
        \ar@<-4pt>[r]
         &
         P \times G
         \ar[d]
         \ar@<+2pt>[r] \ar@<-2pt>[r]
       &
       P
       \ar[d]
       \ar@{->>}[r]
       &
       X
       \ar[d]^g
       \\
       \cdots
       \ar@<+6pt>[r] \ar@<+2pt>[r] \ar@<-2pt>[r] \ar@<-6pt>[r]
       &
       G \times G
       \ar@<+4pt>[r]
       \ar[r]
        \ar@<-4pt>[r]
         &
         G
         \ar@<+2pt>[r] \ar@<-2pt>[r]
       &
       {*}
       \ar@{->>}[r]
       &
       \mathbf{B}G
    }
  $$
  where each square and each composite rectangle is an
  $\infty$-pullback.
  This exhibits the $G$-action on $P$.
  Since $* \to \mathbf{B}G$
  is an effective epimorphism, so is its $\infty$-pullback
  $P \to X$. Since, by the $\infty$-Giraud theorem, $\infty$-colimits are preserved
  by $\infty$-pullbacks we have that $P \to X$ exhibits the
  $\infty$-colimit $X \simeq P/\!/G$.
\endofproof

%
%
 In practice, $\infty$-toposes appear in roughly two guises, analogous to what for ordinary toposes are known as {\it petit toposes} such as sheaves on a fixed scheme or manifold, and {\it gros toposes} such as sheaves on a site of all schemes or all manifolds. Toposes in classical discussions tend to be taken in the  petit form, but often these may be understood as slices (typically {\'e}tale-slices) of gros toposes -- we come back to this in Section~\ref{section.InfinityGerbes} when we connect to the classical discussion of gerbes and higher gerbes.

 One way to formalize the notion of {\it gros $\infty$-toposes} is to characterize their terminal object $\ast$ as being {\it point-like}. There are several ways to make this precise, but in the present context the pertinent condition would be that the map $\ast \to \mathbf{B}G$ in Rem. \ref{PointIntoBGIsEffectiveEpimorphism} is the {\it unique} map of this form, up to homotopy. In this case we have:
%
%

\begin{lemma}\label{TrivialityOfGBundlesOverPoints}
  For $P \to X$ a $G$-principal $\infty$-bundle
  obtained as in Proposition~\ref{BundleStructureOnInfinityFiber},
%
%
in an $\infty$-topos for which there is an essentially unique map $\ast \to \mathbf{B}G$,
%
%
and for
  $x : * \to X$ any point of $X$, we have
  a canonical equivalence
  $$
    \xymatrix{
	    x^* P \ar[r]^{\simeq} & G
	}
  $$
  between the fiber $x^*P$ and the $\infty$-group object $G$.
\end{lemma}
\proof
  This follows from the pasting law for $\infty$-pullbacks, which
  gives the diagram
  $$
    \xymatrix{
      G \ar[d] \ar[r] & P \ar[d]\ar[r] & {*} \ar[d]
      \\
      {*} \ar[r]^x & X \ar[r]^g & \mathbf{B}G
    }
  $$
  in which both squares as well as the total rectangle are
  $\infty$-pullbacks.
%
%
Here the bottom composite is equivalent to the given base point $\ast \to \mathbf{B}G$ by the assumption that all maps of this form are homotopic, and hence the outer rectangle is a pullback by Rem. \ref{PointIntoBGIsEffectiveEpimorphism}.
\endofproof
\begin{definition}
  \label{TrivialGBundle}
  The \emph{trivial} $G$-principal $\infty$-bundle $(P \to X) \simeq (X \times G \to X)$
  is, up to equivalence, the one obtained via Proposition~\ref{BundleStructureOnInfinityFiber}
  from the morphism $X \to * \to \mathbf{B}G$.
\end{definition}
\begin{proposition}
  \label{PullbackOfInfinityTorsors}
  For $P \to X$ a $G$-principal $\infty$-bundle and $Y \to X$ any morphism, the
  $\infty$-pullback $Y \times_X P$ naturally inherits the structure of
  a $G$-principal $\infty$-bundle.
\end{proposition}
\proof
  This uses the same kind of argument as in Proposition~\ref{BundleStructureOnInfinityFiber}.
\endofproof
In fact this is the special case of the pullback of what we will see below in
Proposition~\ref{LocalTrivialityImpliesCocycle} is the
  universal $G$-principal $\infty$-bundle $*\to \mathbf{B}G$.
\begin{proposition}
  \label{EveryGBundleIsLocallyTrivial}
  Every $G$-principal $\infty$-bundle is locally trivial, that is
  there exists an effective epimorphism $\xymatrix{U \ar@{->>}[r] & X}$
  and an equivalence of $G$-principal $\infty$-bundles
  $$
    U \times_X P \simeq U \times G
  $$
  from the pullback of $P$ (Proposition~\ref{PullbackOfInfinityTorsors})
  to the trivial $G$-principal $\infty$-bundle over $U$ (Definition~\ref{TrivialGBundle}).
\end{proposition}
\proof
  For $P \to X$ a $G$-principal $\infty$-bundle, it is, by
  Remark~\ref{GBundlesAreEffectiveEpimorphisms}, itself an effective
  epimorphism. The pullback of the $G$-bundle to its
  own total space along this morphism is trivial,
   by the principality condition
  (Proposition~\ref{PrincipalityCondition}). Hence setting $U := P$
  proves the claim.
\endofproof
\begin{proposition}
  \label{LocalTrivialityImpliesCocycle}
  For every $G$-principal $\infty$-bundle $P \to X$ the square
  $$
    \xymatrix{
       & P \ar[d] \ar[r] & {*} \ar[d]
       \\
       X \ar@{}[r]|<<<\simeq
       & \varinjlim_n (P \times G^{\times_n})
       \ar[r]
       &
       \varinjlim_n G^{\times_n}
       \ar@{}[r]|\simeq
       &
       \mathbf{B}G
    }
  $$
  is an $\infty$-pullback diagram.
\end{proposition}
\proof
  Let $U \to X$ be an effective epimorphism
  such that $P \to X$ pulled back to $U$ becomes the trivial $G$-principal
  $\infty$-bundle. By Proposition~\ref{EveryGBundleIsLocallyTrivial} this exists.
  By definition of morphism of $G$-actions and
  by functoriality of the $\infty$-colimit, this induces
  a morphism in ${\mathbf{H}^{\Delta[1]}}_{/(* \to \mathbf{B}G)}$
  corresponding to the diagram
  $$
    \raisebox{20pt}{
    \xymatrix{
	  U \times G \ar@{->>}[r] \ar@{->>}[d] & P \ar[r] \ar@{->>}[d] & {*} \ar@{->>}[d]^{\mathrm{pt}}
	  \\
	  U \ar@{->>}[r] & X \ar[r] & \mathbf{B}G
	}
	}
	\;\;
	\simeq
	\;\;
    \raisebox{20pt}{
    \xymatrix{
	  U \times G \ar@{->>}[rr] \ar@{->>}[d] & & {*} \ar@{->>}[d]^{\mathrm{pt}}
	  \\
	  U \ar[r] & {*} \ar[r]^{\mathrm{pt}} & \mathbf{B}G
	}
	}
  $$
  in $\mathbf{H}$.
  By assumption, in this diagram the outer rectangles and the square on the very left
  are $\infty$-pullbacks. We need to show that
  the right square on the left is also an $\infty$-pullback.

  Since $U \to X$ is an effective epimorphism by assumption, and since these are
  stable under $\infty$-pullback, $U \times G \to P$ is also an effective epimorphism,
  as indicated. This means that
  $$
    P \simeq {\varinjlim_n}\, (U \times G)^{\times^{n+1}_P}
	\,.
  $$
  We claim that for all $n \in \mathbb{N}$ the fiber products in the colimit on the right
  are naturally equivalent to $(U^{\times^{n+1}_X}) \times G$. For $n = 0$ this is
  clearly true. Assume then by induction that
  it holds for some $n \in \mathbb{N}$. Then
  with the pasting law (Proposition~\ref{PastingLawForPullbacks}) we find  an
  $\infty$-pullback diagram of the form
  $$
    \raisebox{30pt}{
    \xymatrix{
	  (U^{\times^{n+1}_X}) \times G 	
	  \ar@{}[r]|\simeq
	  &
	  (U \times G)^{\times^{n+1}_P}
	  \ar[r] \ar[d]
	  &
	  (U \times G)^{\times^n_P}
      \ar[d]
	  \ar@{}[r]|{\simeq}
	  &
	  (U^{\times^{n}_X}) \times G
	  \\
	  & U \times G \ar[r] \ar[d] & P \ar[d]
	  \\
	  & U \ar[r] & X.
	}
	}
  $$
  This completes the induction.
  With this the above expression for $P$ becomes
  $$
    \begin{aligned}
      P & \simeq {\varinjlim_n}\,  (U^{\times^{n+1}_X}) \times G
	   \\
	    & \simeq {\varinjlim_n}  \,\mathrm{pt}^* \, (U^{\times^{n+1}_X})
	   \\
	    & \simeq \mathrm{pt}^* \, {\varinjlim_n}\,    (U^{\times^{n+1}_X}) 	
	  \\
	    & \simeq \mathrm{pt}^* \, X,
	\end{aligned}
  $$
  where we have used that by the second $\infty$-Giraud axiom (Definition~\ref{GiraudRezkLurieAxioms})
  we may take the $\infty$-pullback out of the $\infty$-colimit and
  where in the last step we
  used again the assumption that $U \to X$ is an effective epimorphism.
\endofproof
\begin{example}
  The fiber sequence
  $$
    \xymatrix{
	  G \ar[r] & {*} \ar[d]
	  \\
	  & \mathbf{B}G
	}
  $$
  which exhibits the delooping $\mathbf{B}G$ of $G$ according to
  Theorem \ref{DeloopingTheorem} is a $G$-principal $\infty$-bundle
  over $\mathbf{B}G$, with \emph{trivial} $G$-action on its total space
  $*$. Proposition \ref{LocalTrivialityImpliesCocycle} says that this is
  the \emph{universal $G$-principal $\infty$-bundle} in that every
  other one arises as an $\infty$-pullback of this one.
  In particular, $\mathbf{B}G$ is a classifying object for $G$-principal
  $\infty$-bundles.

  Below in Theorem \ref{ClassificationOfTwistedGEquivariantBundles}
  this relation is strengthened:
  every \emph{automorphism} of a $G$-principal $\infty$-bundle, and in fact
  its full automorphism $\infty$-group arises from pullback of the above
  universal $G$-principal $\infty$-bundle: $\mathbf{B}G$ is the fine
  \emph{moduli $\infty$-stack} of $G$-principal $\infty$-bundles.

  The traditional definition of universal $G$-principal bundles in terms of
  contractible objects equipped with a free $G$-action has no intrinsic
  meaning in higher topos theory. Instead this appears in
  \emph{presentations} of the general theory in model categories
  (or categories of fibrant objects)
  as \emph{fibrant representatives}
  $\mathbf{E}G \to \mathbf{B}G$ of the above point inclusion.
  This we discuss in \cite{NSSb}.
  \label{UniversalPrincipal}
\end{example}
The main classification Theorem \ref{PrincipalInfinityBundleClassification} below
implies in particular that every morphism in $G\mathrm{Bund}(X)$ is an equivalence.
For emphasis we note how this also follows directly:
\begin{lemma}
  \label{EquivalencesAreDetectedOverEffectiveEpimorphisms}
  Let $\mathbf{H}$ be an $\infty$-topos and let $X$ be an
  object of $\mathbf{H}$.  A morphism
  $f\colon A\to B$ in $\mathbf{H}_{/X}$
  is an equivalence if and only if $p^*f$ is an equivalence in
  $\mathbf{H}_{/Y}$ for any effective epimorphism $p\colon Y\to X$ in
  $\mathbf{H}$.

\end{lemma}
\proof
  It is clear, by functoriality, that $p^* f$ is a weak equivalence if $f$ is.
  Conversely, assume that $p^* f$ is a weak equivalence.
  Since effective epimorphisms as well as
  equivalences are preserved by pullback
  we get a simplicial diagram of the form
  $$
    \raisebox{20pt}{
    \xymatrix{
       \cdots
       \ar@<+4pt>[r]
       \ar@<+0pt>[r]
       \ar@<-4pt>[r]
       &
       p^* A \times_A p^* A
       \ar@<+2pt>[r]
       \ar@<-2pt>[r]
       \ar[d]^\simeq
       &
       p^* A
       \ar[d]^\simeq
       \ar@{->>}[r]
       &
       A
       \ar[d]^f
       \\
       \cdots
       \ar@<+4pt>[r]
       \ar@<+0pt>[r]
       \ar@<-4pt>[r]
       &
       p^* B \times_B p^* B
       \ar@<+2pt>[r]
       \ar@<-2pt>[r]
       &
       p^* B
       \ar@{->>}[r]
       &
       B
    }
	}
  $$
  where the rightmost horizontal morphisms are effective epimorphisms, as indicated.
  By definition of effective epimorphisms this exhibits
  $f$ as an $\infty$-colimit over equivalences, hence as
  an equivalence.
\endofproof
\begin{proposition}
  \label{MorphismsOfInfinityBundlesAreEquivalences}
  Every morphism between $G$-actions over $X$ that are
  $G$-principal $\infty$-bundles over $X$ is an equivalence.
\end{proposition}
\proof
  Since a morphism of $G$-principal bundles
  $P_1 \to P_2$ is a morphism of {\v C}ech nerves that fixes
  their $\infty$-colimit $X$, up to equivalence,
  and since $* \to \mathbf{B}G$ is an effective
  epimorphism,
  we are, by Proposition~\ref{LocalTrivialityImpliesCocycle}, in the situation of
  Lemma~\ref{EquivalencesAreDetectedOverEffectiveEpimorphisms}.
\endofproof
\begin{theorem}
  \label{PrincipalInfinityBundleClassification}
  For all $X, \mathbf{B}G \in \mathbf{H}$ there is a natural
  equivalence of $\infty$-groupoids
  $$
    G \mathrm{Bund}(X)
    \simeq
    \mathbf{H}(X, \mathbf{B}G)
  $$
  which on vertices is the construction of Proposition~\ref{BundleStructureOnInfinityFiber}:
  a bundle $P \to X$ is identified with a morphism
  $X \to \mathbf{B}G$ such that $P \to X \to \mathbf{B}G$ is a fiber
  sequence.
\end{theorem}
We therefore say
\begin{itemize}
  \item $\mathbf{B}G$ is the \emph{classifying object}
  or \emph{moduli $\infty$-stack} for
  $G$-principal $\infty$-bundles;
  \item a morphism $c : X \to \mathbf{B}G$ is a \emph{cocycle}
  for the corresponding $G$-principal $\infty$-bundle and its class
  $[c] \in \mathrm{H}^1(X,G)$ is its
 \emph{characteristic class}.
\end{itemize}
\proof
  By Definitions~\ref{ActionInPrincipal}
  and~\ref{principalbundle} and using
  the refined statement of the third $\infty$-Giraud axiom
  (Theorem~\ref{NaturalThirdGiraud}), the
  $\infty$-groupoid $G\mathrm{Bund}(X)$ of $G$-principal $\infty$-bundles over $X$ is
  equivalent to the fiber over $X$
  of the
  sub-$\infty$-category
  of the slice ${\mathbf{H}^{\Delta[1]}}_{/{(* \to \mathbf{B}G)}}$
  of the arrow $\infty$-topos
  on those squares
  $$
    \xymatrix{
	  P \ar[r] \ar@{->>}[d] & {*} \ar@{->>}[d]
	  \\
	  X \ar[r] & \mathbf{B}G
	}
  $$
  that exhibit $P \to X$ as a $G$-principal $\infty$-bundle.  By
  Proposition~\ref{BundleStructureOnInfinityFiber} and
  Proposition~\ref{LocalTrivialityImpliesCocycle} these are
  $\infty$-pullback squares, hence objects of the full sub-$\infty$-category
  $$
    \mathrm{Cart}({\mathbf{H}^{\Delta[1]}}_{/{(* \to \mathbf{B}G)}})
    \hookrightarrow {\mathbf{H}^{\Delta[1]}}_{/{(* \to \mathbf{B}G)}}
  $$
  of the slice of the arrow category over the point inclusion into $\mathbf{B}G$
  on those morphisms of morphisms (hence squares) which are $\infty$-pullbacks
  (``cartesian'').
  This inclusion is not full, rather the morphisms of $G$-principal $\infty$-bundles
  over $X$ are those morphisms of $G$-actions that fix the base $X$, up to equivalence.
  By the universal property of the homotopy fiber product of $\infty$-categories this means that
  $$
    G \mathrm{Bund}(X) \simeq
	 \mathrm{Cart}({\mathbf{H}^{\Delta[1]}}_{/{(* \to \mathbf{B}G)}}) \times_{\mathbf{H}} \{X\}
	 \,.
  $$
  Now by the universality of the $\infty$-pullback in $\mathbf{H}$
  the morphisms between two Cartesian squares are fixed by their value on the underlying
  co-span $X \to \mathbf{B}G \leftarrow \ast$ and since in the above $\ast \to \mathbf{B}G$
  is held fixed, they are fully determined by their value on $X$,
  so that the above is equivalent to
  $$
  \mathbf{H}_{/\mathbf{B}G} \times_{\mathbf{H}} \{X\}
	\,.
  $$
  Specifically, if one is to do this argument in terms of model categories: choose a model structure for
  $\mathbf{H}$ in which all objects are cofibrant, choose a fibrant representative
  for $\mathbf{B}G$ and a fibration resolution $\mathbf{E}G \to \mathbf{B}G$
  of the universal $G$-bundle. Then the slice model structure of the arrow model structure
  over this presents the slice in question and the statement follows from the analogous
  1-categorical statement.

  This last expression finally is equivalent to
  $$
	\mathbf{H}(X, \mathbf{B}G)
	\,.
  $$
  To see this for instance in terms of quasi-categories: the projection
  $\mathbf{H}_{/\mathbf{B}G} \to \mathbf{H}$ is a fibration by
  Proposition 2.1.2.1 and 4.2.1.6 in \cite{Lurie}, hence the homotopy fiber
  $\mathbf{H}_{/\mathbf{B}G} \times_{\mathbf{X}} \{X\}$
  is the ordinary fiber of quasi-categories. This is manifestly
  the
  $\mathrm{Hom}^R_{\mathbf{H}}(X, \mathbf{B}G)$ from Proposition 1.2.2.3 of \cite{Lurie}.
  Finally, by Proposition 2.2.4.1 there, this is equivalent to $\mathbf{H}(X,\mathbf{B}G)$.
\endofproof
\begin{corollary}
  Equivalence classes of $G$-principal $\infty$-bundles over $X$ are
  in natural bijection with the degree-1 $G$-cohomology of $X$:
  $$
    G \mathrm{Bund}(X)_{/\sim} \simeq H^1(X, G)
	\,.
  $$
\end{corollary}
\proof
  By Definition \ref{cohomology} this is the restriction of
  the equivalence $G \mathrm{Bund}(X) \simeq \mathbf{H}(X, \mathbf{B}G)$ to
  connected components.
\endofproof

\section{Twisted bundles and twisted cohomology}
\label{StrucTwistedCohomology}

We show here how the general notion of cohomology in an
$\infty$-topos, considered above in Section~\ref{StrucCohomology}, subsumes the notion of
\emph{twisted cohomology} and we discuss the corresponding
geometric structures classified by twisted cohomology:
\emph{extensions} of principal $\infty$-bundles and \emph{twisted $\infty$-bundles}.

Whereas ordinary cohomology is given by a derived hom-$\infty$-groupoid,
twisted cohomology is given by the $\infty$-groupoid of
\emph{sections of
a local coefficient bundle} in an $\infty$-topos,
which in turn is
an \emph{associated $\infty$-bundle} induced via a representation
of an $\infty$-group $G$ from a $G$-principal $\infty$-bundle
(this is a geometric and unstable variant of the picture
of twisted cohomology developed in \cite{AndoBlumbergGepner,MaySigurdsson}).

It is fairly immediate that, given a \emph{universal} local coefficient bundle
associated to a universal principal $\infty$-bundle,
the induced twisted cohomology is equivalently ordinary
cohomology in the corresponding slice $\infty$-topos. This
identification provides a clean formulation of the contravariance
of twisted cocycles. However, a universal coefficient bundle
is a pointed connected object in the slice $\infty$-topos only
when it is a trivial bundle, so that twisted cohomology does not classify
principal $\infty$-bundles in the slice. We show below that instead
it classifies \emph{twisted principal $\infty$-bundles}, which are
natural structures that generalize the twisted bundles familiar from
twisted K-theory.
Finally, we observe that twisted cohomology in an $\infty$-topos
equivalently classifies extensions of structure groups
of principal $\infty$-bundles.

A wealth of structures turn out to be special cases of
nonabelian twisted cohomology and of twisted
principal $\infty$-bundles and their study benefits from the
general theory of twisted cohomology.

\subsection{Actions and associated $\infty$-bundles}
\label{StrucRepresentations}

Let $\mathbf{H}$ be an $\infty$-topos, $G \in \mathrm{Grp}(\mathbf{H})$
an $\infty$-group.
Fix an action $\rho : V \times G \to V$  on an object $V\in \mathbf{H}$
as in Definition \ref{ActionInPrincipal}.
We discuss the induced notion of \emph{$\rho$-associated $V$-fiber $\infty$-bundles}.
We show that there is a \emph{universal} $\rho$-associated $V$-fiber bundle over
$\mathbf{B}G$ and observe that under Theorem \ref{PrincipalInfinityBundleClassification}
this is effectively identified with the action itself. Accordingly, we also further discuss
$\infty$-actions as such.

\medskip

\begin{definition}
  For $V,X \in \mathbf{H}$ any two objects,
a \emph{$V$-fiber $\infty$-bundle} over $X$ is a morphism $E \to X$,
such that there is an effective epimorphism
$\xymatrix{U \ar@{->>}[r] & X}$ and an $\infty$-pullback of the form
$$
  \raisebox{20pt}{
  \xymatrix{
    U \times V \ar[r] \ar[d] & E \ar[d]
	\\
	U \ar@{->>}[r] & X\, .
  }
  }
$$
 \label{FiberBundle}
\end{definition}
We say that $E \to X$ locally trivializes with respect to $U$.
As usual, we often say \emph{$V$-bundle} for short.

\begin{definition}
  For $P \to X$ a $G$-principal $\infty$-bundle, we write
  $$
    P \times_G V := (P\times V)/\!/G
  $$
  for the $\infty$-quotient of the diagonal $\infty$-action of $G$ on $P \times V$.
  Equipped with the canonical morphism
  $P \times_G V \to X$ we call this the $\infty$-bundle \emph{ $\rho$-associated} to $P$.
  \label{AssociatedBundle}
\end{definition}
\begin{remark}
  The diagonal $G$-action on $P \times V$ is the product in
  $G \mathrm{Action}(\mathbf{H})$ of the given actions on $P$ and on $V$.
  Since $G\mathrm{Action}(\mathbf{H})$ is a full sub-$\infty$-category of a slice
  category of a functor category, the product is given by a degreewise
  pullback in $\mathbf{H}$:
  $$
    \raisebox{20pt}{
    \xymatrix{
	  P \times V \times G^{\times_n}
	  \ar[r]
	  \ar[d]
	  &
	  V \times G^{\times_n}
	  \ar[d]
	  \\
	  P \times G^{\times_n}
	  \ar[r]
	  &
	  G^{\times_n}\,.
	}
	}
  $$
  and so
  $$
    P \times_G V \simeq \varinjlim_n (P \times V \times G^{\times_n})
	\,.
  $$
  The canonical bundle morphism of the corresponding $\rho$-associated
  $\infty$-bundle  is the realization of the left morphism of this diagram:
  $$
    \raisebox{20pt}{
    \xymatrix{
	  P \times_G V
	  \ar@{}[r]|<<<<{:=}
	  \ar[d]
	  &
	  \varinjlim_n (P \times V \times G^{\times_n})
	  \ar[d]
	  \\
	  X \ar@{}[r]|<<<<<<<<<{\simeq} &
	  \varinjlim_n (P \times G^{\times_n})\,.
	}
	}
  $$
  \label{ProductActionByPullback}
\end{remark}
\begin{example}
By Theorem \ref{PrincipalInfinityBundleClassification} every $\infty$-group action
$\rho : V \times G \to V$ has a classifying morphism $\mathbf{c}$ defined on its homotopy
quotient, which fits into a fiber sequence of the form
$$
  \raisebox{20pt}{
  \xymatrix{
     V \ar[r] & V/\!/G \ar[d]^{\mathbf{c}}
     \\	
	 & \mathbf{B}G\,.
  }
  }
$$

  Regarded as an
  $\infty$-bundle, this is
  $\rho$-associated to the universal $G$-principal $\infty$-bundle
  $\xymatrix{{*} \ar[r] & \mathbf{B}G}$ from Example \ref{UniversalPrincipal}:
  $$
    V/\!/G \simeq {*} \times_G V
	\,.
  $$
  \label{ActionGroupoidIsRhoAssociated}
\end{example}
\begin{lemma}
  The realization functor $\varinjlim : \mathrm{Grpd}(\mathbf{H}) \to \mathbf{H}$
  preserves the $\infty$-pullback of Remark \ref{ProductActionByPullback}:
  $$
    P \times_G V \simeq \varinjlim_n (P \times V \times G^{\times_n})
	\simeq
	(\varinjlim_n P \times G^{\times_n}) \times_{(\varinjlim_n G^{\times_n})} (\varinjlim_n V \times G^{\times_n})
	\,.
  $$
  \label{RealizationPreservesProductOfRepresentations}
\end{lemma}
\proof
  Generally, let $X \to Y \leftarrow Z \in \mathrm{Grpd}(\mathbf{H})$ be a
  diagram of groupoid objects, such that in the induced diagram
  $$
    \xymatrix{
	  X_0 \ar[r] \ar@{->>}[d] & Y_0 \ar@{<-}[r] \ar@{->>}[d] & Z_0 \ar@{->>}[d]
	  \\
	  \varinjlim_n X_n \ar[r] & \varinjlim_n Y_n \ar@{<-}[r] & \varinjlim_n Z_n
	}
  $$
  the left square is an $\infty$-pullback. By the third
  $\infty$-Giraud axiom (Definition~\ref{GiraudRezkLurieAxioms}) the vertical
  morphisms are effective epis, as indicated.
  By assumption we have a pasting of $\infty$-pullbacks as shown on the
  left of the following diagram, and by
  the pasting law (Proposition \ref{PastingLawForPullbacks}) this is equivalent to
  the pasting shown on the right:
  $$
    \raisebox{38pt}{
    \xymatrix{
	  X_0 \times_{Y_0} Z_0 \ar[r] \ar[d] & Z_0 \ar[d]
	  \\
	  X_0 \ar[r] \ar[d] & Y_0 \ar[d]
	  \\
	  \varinjlim_n X_n \ar[r] & \varinjlim_n Y_n
	}
	}
	\;\;\;
	\simeq
	\;\;\;
    \raisebox{38pt}{
    \xymatrix{
	  X_0 \times_{Y_0} Z_0 \ar[r] \ar@{->>}[d] & Z_0 \ar@{->>}[d]
	  \\
	  (\varinjlim_n X_n) \times_{(\varinjlim_n Y_n)} (\varinjlim_n Z_n) \ar[r] \ar[d] &
	  \varinjlim_n Z_n \ar[d]
	  \\
	  \varinjlim_n X_n \ar[r] & \varinjlim_n Y_n.
	}
	}
  $$
Since effective epimorphisms are stable under $\infty$-pullback, this identifies
the canonical morphism
$$
  X_0 \times_{Y_0} Z_0
  \to
  (\varinjlim_n X_n) \times_{(\varinjlim_n Y_n)} (\varinjlim_n Z_n)
$$
as an effective epimorphism, as indicated.

Since $\infty$-limits commute over each other, the {\v C}ech nerve of this morphism
is the groupoid object $[n] \mapsto X_n \times_{Y_n} Z_n$.
Therefore the third $\infty$-Giraud axiom now says that $\varinjlim$ preserves the
$\infty$-pullback of groupoid objects:
$$
  \varinjlim (X \times_Y Z)
   \simeq
  \varinjlim_n (X_n \times_{Y_n} Z_n )
   \simeq
  (\varinjlim_n X_n) \times_{(\varinjlim_n Y_n)} (\varinjlim_n Z_n)
  \,.
$$

Consider this now in the special case that $X \to Y \leftarrow Z$ is
$(P \times G^{\times_\bullet}) \to G^{\times_\bullet} \leftarrow (V \times G^{\times_\bullet})$.
Theorem \ref{PrincipalInfinityBundleClassification} implies that the initial assumption above is
met, in that $P \simeq (P/\!/G) \times_{*/\!/G} {*} \simeq X \times_{\mathbf{B}G} {*}$,
and so the claim follows.
\endofproof
\begin{proposition}
  For $g_X : X \to \mathbf{B}G$ a morphism and $P \to X$
  the corresponding $G$-principal $\infty$-bundle according to Theorem
  \ref{PrincipalInfinityBundleClassification},
  there is a natural equivalence
  $$
    g_X^*(V/\!/G) \simeq P \times_G V
  $$
  over $X$, between the pullback of the
  $\rho$-associated $\infty$-bundle
  $\xymatrix{V/\!/G \ar[r]^{\mathbf{c}} & \mathbf{B}G}$
  of Example \ref{ActionGroupoidIsRhoAssociated}
  and the $\infty$-bundle $\rho$-associated to $P$ by Definition \ref{AssociatedBundle}.
  \label{UniversalAssociatedBundle}
\end{proposition}
\proof
 By Remark \ref{ProductActionByPullback} the product action is given by the
 pullback
 $$
   \xymatrix{
     P \times V \times G^{\times_\bullet}
	 \ar[r]
	 \ar[d]
	 &
	 V \times G^{\times_\bullet}
	 \ar[d]
	 \\
	 P \times G^{\times_\bullet} \ar[r] & G^{\times_\bullet}
   }
 $$
 in $\mathbf{H}^{\Delta^{\mathrm{op}}}$.
 By Lemma $\ref{RealizationPreservesProductOfRepresentations}$ the realization functor
 preserves this $\infty$-pullback. By
 Remark \ref{ProductActionByPullback} it sends the left morphism to the
 associated bundle, and by Theorem \ref{PrincipalInfinityBundleClassification}
 it sends the bottom morphism to $g_X$. Therefore it produces an $\infty$-pullback
 diagram of the form
 $$
   \raisebox{20pt}{
   \xymatrix{
     V \times_G P \ar[r] \ar[d] & V/\!/G \ar[d]^{\mathbf{c}}
	 \\
	 X \ar[r]^{g_X} & \mathbf{B}G\,.
   }
   }
 $$
\endofproof
\begin{remark}
  This says that $\xymatrix{V/\!/G \ar[r]^{\mathbf{c}} & \mathbf{B}G}$ is both,
  the $V$-fiber $\infty$-bundle
  $\rho$-associated to the universal $G$-principal $\infty$-bundle, Example
  \ref{ActionGroupoidIsRhoAssociated},
  as well as the universal $\infty$-bundle for $\rho$-associated $\infty$-bundles.
  \label{RhoAssociatedToUniversalIsUniversalVBundle}
\end{remark}
\begin{proposition}
  Every $\rho$-associated $\infty$-bundle is a $V$-fiber $\infty$-bundle,
  Definition \ref{FiberBundle}.
  \label{AssociatedIsFiberBundle}
\end{proposition}
\proof
  Let $P \times_G V \to X$ be a $\rho$-associated $\infty$-bundle.
  By the previous Proposition \ref{UniversalAssociatedBundle} it is
  the pullback $g_X^* (V/\!/G)$ of the universal $\rho$-associated bundle.
  By Proposition \ref{EveryGBundleIsLocallyTrivial} there exists an
  effective epimorphism $\xymatrix{U \ar@{->>}[r] & X}$ over which
  $P$ trivializes, hence such that $g_X|_U$ factors through the point, up
  to equivalence. In summary and by the pasting law, Proposition \ref{PastingLawForPullbacks},
  this gives a pasting of $\infty$-pullbacks of the form
  $$
    \raisebox{20pt}{
    \xymatrix@R=8pt{
	  U \times V
	  \ar[dd]
	  \ar[r]
	  &
	  P \times_G V
	  \ar[r]
	  \ar[dd]
	  &
	  V/\!/G
	  \ar[dd]
	  \\
	  \\
	  U \ar@{->>}[r]
	  \ar[dr]
	  &
	  X
	  \ar[r]^{g_X}
	  &
	  \mathbf{B}G
	  \\
	  & {*}
	  \ar[ur]
	}
	}
  $$
  which exhibits $P \times_G V \to X$ as a $V$-fiber bundle by a local trivialization
  over $U$.
\endofproof

So far this shows that every $\rho$-associated $\infty$-bundle is a
$V$-fiber bundle. We want to show that, conversely, every $V$-fiber bundle
is associated to a principal $\infty$-bundle.
\begin{definition}
  Let $V \in \mathbf{H}$ be a $\kappa$-compact object, for some regular cardinal $\kappa$.
  By the characterization of Proposition \ref{RezkCharacterization}, there exists
  an $\infty$-pullback square in $\mathbf{H}$ of the form
  $$
    \raisebox{20pt}{
    \xymatrix{
	  V \ar[r] \ar[d] & \widehat {\mathrm{Obj}}_\kappa \ar[d]
	  \\
	  {*} \ar[r]^<<<<<<{\vdash V} & \mathrm{Obj}_\kappa
	}
	}
  $$
  Write
  $$
    \mathbf{B}\mathbf{Aut}(V) := \mathrm{im}(\vdash V)
  $$
  for the $\infty$-image, Definition \ref{image},
  of the classifying morphism $\vdash V$ of $V$.
  By definition this comes with an effective epimorphism
  $$
    \xymatrix{
	  {*} \ar@{->>}[r] & \mathbf{B}\mathbf{Aut}(V)
	  \ar@{^{(}->}[r] & \mathrm{Obj}_\kappa
	}
	\,,
  $$
  and hence, by Proposition \ref{InfinityGroupObjectsAsGroupoidObjects},
  it is the delooping of an $\infty$-group
  $$
    \mathbf{Aut}(V) \in \mathrm{Grp}(\mathbf{H})
  $$
  as indicated. We call this the \emph{internal automorphism $\infty$-group} of $V$.

  By the pasting law, Proposition \ref{PastingLawForPullbacks},
  the image factorization gives a pasting
  of $\infty$-pullback diagrams of the form
  $$
    \xymatrix{
	  V \ar[r] \ar[d] & V/\!/\mathbf{Aut}(V) \ar[r] \ar[d]^{\mathbf{c}_V} &
	  \widehat {\mathrm{Obj}}_\kappa \ar[d]
	  \\
	  {*} \ar@{->>}[r]^<<<<<<{\vdash V} & \mathbf{B}\mathbf{Aut}(V) \ar@{^{(}->}[r] &
	  \mathrm{Obj}_\kappa
	}
  $$
  By Theorem~\ref{PrincipalInfinityBundleClassification} this defines a canonical
  $\infty$-action
  $$
    \rho_{\mathbf{Aut}(V)} : V \times \mathbf{Aut}(V) \to V
  $$
  of $\mathbf{Aut}(V)$ on $V$ with homotopy quotient $V/\!/\mathbf{Aut}(V)$
  as indicated.
  \label{InternalAutomorphismGroup}
\end{definition}
\begin{proposition}
  Every $V$-fiber $\infty$-bundle is $\rho_{\mathbf{Aut}(V)}$-associated to an
  $\mathbf{Aut}(V)$-principal $\infty$-bundle.
  \label{VBundleIsAssociated}
\end{proposition}
\proof
  Let $E \to V$ be a $V$-fiber $\infty$-bundle.
  By Definition \ref{FiberBundle} there exists an effective epimorphism
  $\xymatrix{U \ar@{->>}[r] & X}$  along which the bundle trivializes locally.
  It follows
  by the second Axiom in Proposition \ref{RezkCharacterization}
  that on $U$ the morphism $\xymatrix{X \ar[r]^<<<<<{\vdash E} & \mathrm{Obj}_\kappa}$
  which classifies $E \to X$ factors through the point
  $$
    \raisebox{20pt}{
    \xymatrix@R=8pt{
	  U \times V \ar[r]\ar[dd] & E \ar [r] \ar[dd] & \widehat{\mathrm{Obj}}_\kappa \ar[dd]
	  \\
	  \\
	  U \ar@{->>}[r] \ar[dr] & X \ar[r]^<<<<<{\vdash E} & \mathrm{Obj}_\kappa.
	  \\
	  & {*} \ar[ur]_<<<<<{\vdash V}
	}
	}
  $$
  Since the point inclusion, in turn, factors through its $\infty$-image
  $\mathbf{B}\mathbf{Aut}(V)$, Definition \ref{InternalAutomorphismGroup},
  this yields the outer commuting diagram of the following form
  $$
    \raisebox{20pt}{
    \xymatrix{
	  U \ar[r] \ar@{->>}[d] & {*} \ar[r] &  \mathbf{B}\mathbf{Aut}(V) \ar@{^{(}->}[d]
	  \\
	  X \ar[rr]_{\vdash E}
	   \ar@{-->}[urr]^{g}
	  && \mathrm{Obj}_\kappa
	}}
  $$
  By the epi/mono factorization system of Proposition \ref{EpiMonoFactorizationSystem}
  there is a diagonal lift $g$ as indicated. Using again the
  pasting law and by Definition \ref{InternalAutomorphismGroup}
  (and the discussion following that)
  this factorization induces a pasting of $\infty$-pullbacks of the form
  $$
    \raisebox{20pt}{
    \xymatrix{
	  E \ar[r] \ar[d] & V/\!/\mathbf{Aut}(V) \ar[r] \ar[d]^{\mathbf{c}_V} &
	  \widehat {\mathrm{Obj}}_\kappa \ar[d]
	  \\
	  X \ar[r]^<<<<<g & \mathbf{B}\mathbf{Aut}(V) \ar@{^{(}->}[r] & \mathrm{Obj}_\kappa
	}
	}
  $$
  Finally, by Proposition \ref{UniversalAssociatedBundle}, this
  exhibits $E \to X$ as being $\rho_{\mathbf{Aut}(V)}$-associated to the
  $\mathbf{Aut}(V)$-principal $\infty$-bundle with class $[g] \in H^1(X,G)$.
\endofproof
\begin{theorem}
  $V$-fiber $\infty$-bundles over $X \in \mathbf{H}$ are classified by
  $H^1(X, \mathbf{Aut}(V))$.
  \label{VBundleClassification}
\end{theorem}
  Under this classification, the $V$-fiber $\infty$-bundle corresponding to
  $[g] \in H^1(X, \mathbf{Aut}(V))$ is identified,
  up to equivalence, with the $\rho_{\mathbf{Aut}(V)}$-associated $\infty$-bundle
  (as in Definition~\ref{AssociatedBundle}) to the $\mathbf{Aut}(V)$-principal $\infty$-bundle
  corresponding to $[g]$ by Theorem~\ref{PrincipalInfinityBundleClassification}.
 \\
\proof
  By Proposition~\ref{VBundleIsAssociated} every morphism
  $\xymatrix{
    X \ar[r]^<<<<<{\vdash E} & \mathrm{Obj}_\kappa
  }$
  that classifies a small $\infty$-bundle $E \to X$
  which happens to be a $V$-fiber $\infty$-bundle factors via some $g$
  through
  the moduli $\infty$-stack $\mathbf{B}\mathbf{Aut}(V)$ for $\mathbf{Aut}(V)$-principal $\infty$-bundles
  $$
    \xymatrix{
	  X
	  \ar[r]^<<<<<{g}
	  \ar@/_1pc/[rr]_{\vdash E}
	  &
      \mathbf{B}\mathbf{Aut}(V)
	  \ar@{^{(}->}[r]
	  & \mathrm{Obj}_\kappa
     }
	 \,.
   $$
   Therefore it only remains to show that also every homotopy
   $(\vdash E_1) \Rightarrow (\vdash E_2)$ factors
   through a homotopy $g_1 \Rightarrow g_2$.
   This follows by applying the epi/mono lifting property of
   Proposition \ref{EpiMonoFactorizationSystem} to the diagram
   $$
     \xymatrix{
	  X \coprod X \ar[r]^<<<<<{(g_1, g_2)} \ar@{->>}[d] & \mathbf{B}\mathbf{Aut}(V)
	  \ar@{^{(}->}[d]
	  \\
	  X \ar[r] \ar@{-->}[ur] & \mathrm{Obj}_\kappa
	 }
   $$
   The outer diagram exhibits the original homotopy. The left morphism is
   an effective epi (for instance immediately by Proposition \ref{EffectiveEpiIsEpiOn0Truncation}),
   the right morphism is a monomorphism by construction. Therefore the dashed
   lift exists as indicated and so the top left triangular diagram exhibits
   the desired factorizing homotopy.
\endofproof
\begin{remark}
  In the special case that $\mathbf{H} = \mathrm{Grpd}_{\infty}$, the classification
  Theorem \ref{VBundleClassification}
  is classical \cite{Stasheff,May}, traditionally
  stated in (what in modern terminology is) the
  presentation of $\mathrm{Grpd}_{\infty}$ by simplicial sets
  or by topological spaces. Recent discussions include \cite{BlomgrenChacholski}.
  For $\mathbf{H}$ a general 1-localic $\infty$-topos (meaning: with a 1-site of definition),
  the statement of Theorem \ref{VBundleClassification} appears in \cite{Wendt},
  formulated there in terms of the presentation of $\mathbf{H}$ by simplicial presheaves.
  (We discuss the relation of these presentations to the above general abstract result
  in \cite{NSSb}.)
  Finally, one finds that the classification of \emph{$G$-gerbes} \cite{Giraud} and
  \emph{$G$-2-gerbes} in \cite{Breen} is the special case of the general statement,
  for $V = \mathbf{B}G$ and $G$ a 1-truncated $\infty$-group.
  This we discuss below in Section~\ref{StrucInftyGerbes}.
  \label{ReferencesOnClassificationOfVBundles}
\end{remark}

We close this section with a list of some fundamental classes
of examples of $\infty$-actions, or equivalently,
by Remark \ref{RhoAssociatedToUniversalIsUniversalVBundle},
of universal associated $\infty$-bundles.
For doing so we use again that,
by Theorem \ref{PrincipalInfinityBundleClassification}, to give an
$\infty$-action of $G$ on $V$ is equivalent to giving a fiber
sequence of the form $V \to V/\!/G \to \mathbf{B}G$.
Therefore the following list mainly serves to associate a traditional
\emph{name} with a given $\infty$-action.
\begin{example}
    \label{ExamplesOfActions}
  The following are $\infty$-actions.
  \begin{enumerate}
    \item
	  For every $V \in \mathbf{H}$, the fiber sequence
	  $$
	   \raisebox{20pt}{
	   \xymatrix{
	      V
		  \ar[d]^{(\mathrm{id}_V, \mathrm{pt}_{\mathbf{B}G})}
		  \\
		  V \times \mathbf{B}G
		  \ar[r]^-{p_2}
		  &
		  \mathbf{B}G
		}
		}
	  $$
	  is the \emph{trivial $\infty$-action} of $G$ on $V$.
    \item
	 For every $G \in \mathrm{Grp}(\mathbf{H})$, the fiber sequence
	 $$
	   \raisebox{20pt}{
	   \xymatrix{
	     G \ar[d]
		 \\
		 {*} \ar[r] & \mathbf{B}G
	   }
	   }
	 $$
	 which defines $\mathbf{B}G$ by Theorem \ref{DeloopingTheorem}
	 induces the \emph{right action of $G$ on itself}
	 $$
	   * \simeq G/\!/G
	   \,.
	 $$
	 At the same time
	 this sequence, but now regarded as a bundle over $\mathbf{B}G$,
	 is the universal $G$-principal $\infty$-bundle, Remark \ref{UniversalPrincipal}.
   \item
     For every object $X \in \mathbf{H}$ write
	 $$
	   \mathbf{L}X := X \times_{X \times X} X
	 $$
	 for its \emph{free loop space} object, the $\infty$-fiber product of the
	 diagonal on $X$ along itself
	 $$
	   \raisebox{20pt}{
	   \xymatrix{
	     \mathbf{L}X \ar[r] \ar[d]_{\mathrm{ev}_{*}} & X \ar[d]
		 \\
		 X \ar[r] & X \times X
	   }
	   }
	 $$
     For every $G \in \mathrm{Grp}(\mathbf{H})$ there is a fiber sequence
	 $$
	   \raisebox{20pt}{
	   \xymatrix{
	     G
		 \ar[d]
		 \\
		 \mathbf{L}\mathbf{B}G
		 \ar[r]^{\mathrm{ev}_{*}}
		 &
		 \mathbf{B}G
	   }
	   }
	 $$
	 This exhibits the \emph{adjoint action of $G$ on itself}
	 $$
	   \mathbf{L}\mathbf{B}G \simeq G/\!/_{\mathrm{ad}} G
	   \,.
	 $$
	\item
	  For every $V \in \mathbf{H}$ there is the canonical
	  $\infty$-action of the \emph{automorphism $\infty$-group}
	  $$
	    \raisebox{20pt}{
	    \xymatrix{
		  V
		  \ar[d]
		  \\
		  V/\!/\mathbf{Aut}(V)
		  \ar[r]
		  &
		  \mathbf{B}\mathbf{Aut}(V)
		}
		}
	  $$
	  introduced in Definition \ref{InternalAutomorphismGroup}, this exhibits the
	  \emph{automorphism action}.
	\item
	  For $\rho_1, \rho_2 \in \mathbf{H}_{/\mathbf{B}G}$ two
	  $G$-$\infty$-actions on objects $V_1, V_2 \in \mathbf{H}$, respectively,
	  their internal hom $[\rho_1, \rho_2] \in \mathbf{H}_{/\mathbf{B}G}$ in the slice over $\mathbf{B}G$ is a $G$-$\infty$-action on the internal hom $[V_1, V_2] \in \mathbf{H}$:
	  $$
	    \xymatrix{
		  [V_1, V_2]
		  \ar[d]
		  \\
		  [V_1, V_2]/\!/ G
		  \ar[r] & \mathbf{B}G
		}
	  $$
	  hence $[V_1, V_2]/\!/G \simeq \sum_{\mathbf{B}G}[\rho_1, \rho_2]$,
     where $\sum_{\mathbf{B}G} : \mathbf{H}_{/\mathbf{B}G} \to \mathbf{H}$ is the left adjoint
     to pullback along the terminal map.
	  (This follows by the fact that the inverse image of base change along
	  $\mathrm{pt}_{\mathbf{B}G} : * \to \mathbf{B}G$ is a cartesian closed $\infty$-functor
	  and hence preserves internal homs.\footnote{U.S. thanks Mike Shulman for discussion of this point.})
	  This is the \emph{conjugation $\infty$-action} of $G$ on morphisms $V_1 \to V_2$ by
	  pre- and postcomposition with the action of $G$ on $V_1$ and $V_2$, respectively.
  \end{enumerate}
\end{example}

\subsection{Sections and twisted cohomology}
\label{TwistedCohomology}

We discuss a general notion of \emph{twisted cohomology} or
\emph{cohomology with local coefficients} in any $\infty$-topos $\mathbf{H}$,
where the \emph{local coefficient $\infty$-bundles} are
associated $\infty$-bundles as discussed above, and where the cocycles are
\emph{sections} of these local coefficient bundles.

\medskip

\begin{definition}
  Let $p : E \to X$ be any morphism in $\mathbf{H}$, to be regarded as an
  $\infty$-bundle over $X$. A \emph{section} of $E$ is a diagram
  $$
    \raisebox{20pt}{
    \xymatrix@!C=40pt@!R=30pt{
	   & E \ar[d]^p
	   \\
	  X \ar[r]_{\mathrm{id}}^{\ }="t"
	   \ar[ur]^{\sigma}_{\ }="s"
	  & X
	  \ar@{=>}^{\simeq} "s"; "t"
	}
	}
  $$
  (where for emphasis we display the presence of the homotopy filling the diagram).
  The \emph{$\infty$-groupoid of sections} of $E \stackrel{p}{\to} X$ is the homotopy fiber
  $$
    \Gamma_X(E) := \mathbf{H}(X,E) \times_{\mathbf{H}(X,X)} \{\mathrm{id}_X\}
  $$
  of the space of all morphisms $X \to E$ on those that cover the identity on $X$.
  \label{Sections}
\end{definition}
We record two elementary but important lemmas about spaces of sections.
\begin{lemma}
  There is a canonical identification
  $$
    \Gamma_X(E) \simeq \mathbf{H}_{/X}(\mathrm{id}_X, p)
  $$
  of the space of sections of $E \to X$ with the hom-$\infty$-groupoid in the
  slice $\infty$-topos $\mathbf{H}_{/X}$ between the identity on $X$ and the bundle map $p$.
  \label{SectionBySliceMaps}
\end{lemma}
\proof
  For instance by Proposition 5.5.5.12 in \cite{Lurie}.
\endofproof
\begin{lemma}
  Let
  $$
    \xymatrix{
	  E_1 \ar[r] \ar[d]^{p_1} & E_2 \ar[d]^{p_2}
	  \\
	  B_1 \ar[r]^{f} & B_2
	}
  $$
  be an $\infty$-pullback diagram in $\mathbf{H}$
  and let $\xymatrix{X \ar[r]^{g_X} & B_1}$ be any morphism. Then
  post-composition with $f$ induces
  a natural equivalence of hom-$\infty$-groupoids
  $$
    \mathbf{H}_{/B_1}(g_X, p_1)
	\simeq
	\mathbf{H}_{/B_2}(f\circ g_X, p_2)
	\,.
  $$
  \label{SliceMapsIntoPullbacks}
\end{lemma}
\proof
  By Proposition 5.5.5.12 in \cite{Lurie}, the left hand side is given by the homotopy pullback
  $$
    \raisebox{20pt}{
    \xymatrix{
	  \mathbf{H}_{/B_1}(g_X, p_1) \ar[r] \ar[d] & \mathbf{H}(X, E_1)
	    \ar[d]^{\mathbf{H}(X,p_1)}
	  \\
	  \{g_X\} \ar[r] & \mathbf{H}(X,B_1)\,.
	}
	}
  $$
  Since the hom-$\infty$-functor
  $\mathbf{H}(X,-) : \mathbf{H} \to \mathrm{Grpd}_{\infty}$
  preserves the $\infty$-pullback $E_1 \simeq f^* E_2$, this
  extends to a pasting of $\infty$-pullbacks, which by the pasting law
  (Proposition \ref{PastingLawForPullbacks}) is
  $$
    \raisebox{20pt}{
    \xymatrix{
	  \mathbf{H}_{/B_1}(g_X, p_1) \ar[r] \ar[d]
	     &
	  \mathbf{H}(X, E_1) \ar[d]^{\mathbf{H}(X,p_1)} \ar[r]
	  &
	  \mathbf{H}(X, E_2)
	  \ar[d]^{\mathbf{H}(X, p_2)}
	  \\
	  \{g_X\} \ar[r] & \mathbf{H}(X,B_1)
	  \ar[r]_-{\mathbf{H}(X, f)} &
      \mathbf{H}(X, B_2)	
	}
	}
	\;\;\;
	\simeq
	\;\;\;
	\raisebox{20pt}{
	\xymatrix{
	  \mathbf{H}_{/B_2}(f \circ g_X, p_2)
	  \ar[r]
	  \ar[d]
	  &
	  \mathbf{H}(X, E_2)
	  \ar[d]^{\mathbf{H}(X, p_2)}
	  \\
	  \{f\circ g_X\}
	  \ar[r]
	  &
	  \mathbf{H}(X, B_2).
	}
	}
  $$
\endofproof
Fix now an $\infty$-group $G \in \mathrm{Grp}(\mathbf{H})$ and an
$\infty$-action $\rho : V \times G \to V$. Write
$$
  \xymatrix{
    V \ar[r] & V/\!/G \ar[d]^{\mathbf{c}}
	\\
	& \mathbf{B}G
  }
$$
for the corresponding \emph{universal $\rho$-associated $\infty$-bundle}
as discussed in Section~\ref{StrucRepresentations}.
\begin{proposition}
  For $g_X : X \to \mathbf{B}G$ a cocycle and $P \to X$ the corresponding
$G$-principal $\infty$-bundle according to Theorem \ref{PrincipalInfinityBundleClassification},
there is a natural equivalence
$$
  \Gamma_X(P \times_G V) \simeq \mathbf{H}_{/\mathbf{B}G}(g_X, \mathbf{c})
$$
between the space of sections of the corresponding $\rho$-associated $V$-bundle
(as in Definition~\ref{AssociatedBundle}) and the hom-$\infty$-groupoid of the slice
$\infty$-topos of $\mathbf{H}$ over $\mathbf{B}G$, between $g_X$ and $\mathbf{c}$.
Schematically:
\[
\left\{
	\begin{xy}
		(10,10)*+{E}="1";
		(-10,-10)*+{X}="2";
		(10,-10)*+{X}="3";
		(5,0)*+{}="4";
		(3,-5)*+{}="5";
		{\ar^-{p} "1";"3"};
		{\ar^-{\sigma} "2";"1"};
		{\ar_-{\mathrm{id}} "2";"3"};	
		{\ar@{=>}^-{\simeq} "4";"5"};
	\end{xy}
\right\}
 \;\;
    \simeq
  \;\;
\left\{
\begin{xy}
		(10,10)*+{V/\!/G}="1";
		(-10,-10)*+{X}="2";
		(10,-10)*+{\mathbf{B}G}="3";
		(5,0)*+{}="4";
		(3,-5)*+{}="5";
		{\ar^-{\mathbf{c}} "1";"3"};
		{\ar^-{\sigma} "2";"1"};
		{\ar_-{\mathrm{g_X}} "2";"3"};	
		{\ar@{=>}^-{\simeq} "4";"5"};
	\end{xy}
\right\}
\]
\label{SectionsAndSliceHoms}
\end{proposition}
\proof
  By Lemma \ref{SectionBySliceMaps} and Lemma \ref{SliceMapsIntoPullbacks}.
\endofproof
\begin{corollary}
  The $\infty$-groupoid of sections of the associated bundle $P \times_G V$
  is naturally equivalent to the $\infty$-groupoid of morphisms of
  $G$-actions $P \to V$:
  $$
    \Gamma_X(P \times_G V)
	\simeq
	G\mathrm{Action}(\mathbf{H})(P,V)
	\,.
  $$
  \label{SectionsOfAssociatedIsActionMorphism}
\end{corollary}
\proof
  Using Proposition \ref{SectionsAndSliceHoms} with
  Theorem \ref{PrincipalInfinityBundleClassification} and with the definition of $\mathbf{c}$
  shows that a section $\sigma \in \Gamma_X(P \times_G V)$ is equivalently a
  morphism of $G$-actions
  $\bar \sigma : P \to V$ from the total space of the $G$-principal bundle,
  namely the morphism
  on homotopy fibers induced from the commuting square underlying $\sigma$ when regarded
  as an element of $\mathbf{H}_{/\mathbf{B}G}(g_X, \mathbf{p})$:
  $$
    \raisebox{30pt}{
    \xymatrix{
	  P \ar[d] \ar[r]^{\bar \sigma} & V \ar[d]
	  \\
	  X \ar[d]^{g_X}\ar[r]^\sigma & V/\!/G \ar[d]^{\mathbf{c}}
	  \\
	  \mathbf{B}G
	  \ar[r]^=
	  &
	  \mathbf{B}G
	}
	}
	\,.
  $$
  By Definition \ref{ActionInPrincipal} of $G \mathrm{Action}(\mathbf{H})$
  as the full sub-$\infty$-category of the slice
  $\mathrm{Grpd}(\mathbf{H})_{/*/\!/G}$, this establishes the equivalence.
\endofproof
\begin{proposition}
  If in the above the cocycle $g_X$ is trivializable, in the sense that it factors through the
  point $* \to \mathbf{B}G$ (equivalently if its class $[g_X] \in H^1(X,G)$ is
  trivial) then there is an equivalence
  $$
    \mathbf{H}_{/\mathbf{B}G}(g_X, \mathbf{c})
	\simeq
	\mathbf{H}(X,V)
	\,.
  $$
  \label{TwistedCohomologyIsLocallyVCohomology}
\end{proposition}
\proof
  In this case the homotopy pullback on the right in the proof of
  Proposition \ref{SectionsAndSliceHoms} is
  $$
	\raisebox{20pt}{
	\xymatrix{
	  \mathbf{H}_{/\mathbf{B}G}(g_X, \mathbf{c})
	  \ar@{}[r]|{\simeq}
	  &
	  \mathbf{H}(X, V)
	  \ar[r]
	  \ar[d]
	  &
	  \mathbf{H}(X, V/\!/G)
	  \ar[d]^{\mathbf{H}(X, \mathbf{c})}
	  \\
	  \{g_X\}
	  \ar@{}[r]|\simeq
	  &
	  \mathbf{H}(X, {*})
	  \ar[r]
	  &
	  \mathbf{H}(X, \mathbf{B}G)
	}
	}
  $$
  using that $V \to V/\!/G \stackrel{\mathbf{c}}{\to} \mathbf{B}G$ is a fiber
  sequence by definition, and that $\mathbf{H}(X,-)$ preserves this fiber sequence.
\endofproof
\begin{remark}
  Since by Proposition~\ref{EveryGBundleIsLocallyTrivial}
  every cocycle $g_X$ trivializes locally over some cover
  $\xymatrix{U \ar@{->>}[r] & X}$
  and equivalently, by Proposition \ref{AssociatedIsFiberBundle}, every
  $\infty$-bundle $P \times_G V$ trivializes locally,
  Proposition \ref{TwistedCohomologyIsLocallyVCohomology} says that
  elements $ \sigma \in \Gamma_X(P \times_G V) \simeq \mathbf{H}_{/\mathbf{B}G}(g_X, \mathbf{c})$
  \emph{locally} are morphisms $\sigma|_U : U \to V$ with values in $V$.
   They fail to be so \emph{globally} to the extent that $[g_X] \in H^1(X, G)$ is non-trivial,
   hence to the extent that $P \times_G V \to X$ is non-trivial.
\end{remark}
This motivates the following definition.
\begin{definition}
  We say that
  the $\infty$-groupoid
  $\Gamma_X(P \times_G V) \simeq \mathbf{H}_{/\mathbf{B}G}(g_X, \mathbf{c})$
  from Proposition \ref{SectionsAndSliceHoms}
  is the $\infty$-groupoid of \emph{$[g_X]$-twisted cocycles} with values in $V$, with respect to the
  \emph{local coefficient $\infty$-bundle} $V/\!/G \stackrel{\mathbf{c}}{\to} \mathbf{B}G$.

  Accordingly, its set of connected components we call the
  \emph{$[g_X]$-twisted $V$-cohomology} with respect to the local coefficient bundle $\mathbf{c}$
  and write:
  $$
    H^{[g_X]}(X, V) := \pi_0 \mathbf{H}_{/\mathbf{B}G}(g_X, \mathbf{c})
	\,.
  $$
  \label{TwistedCohomologyInOvertopos}
\end{definition}
\begin{remark}
  The perspective that twisted cohomology is the theory of sections of
  associated bundles whose fibers are classifying spaces is maybe most famous
  for the case of twisted K-theory, where it was described in this form in
  \cite{Rosenberg}.
  But already the old theory of
  \emph{ordinary cohomology with local coefficients} is of this form,
  as is made manifest in \cite{BFG}
  .

  A proposal for a comprehensive theory in terms of bundles of topological spaces is
  in \cite{MaySigurdsson} and a systematic formulation
  in $\infty$-category theory and
  for the case of multiplicative generalized cohomology theories is in
  \cite{AndoBlumbergGepner}. The formulation above refines this, unstably,
  to geometric cohomology
  theories/(nonabelian) sheaf hypercohomology,
  hence from bundles of classifying spaces to
  $\infty$-bundles of moduli $\infty$-stacks.

  A wealth of examples and applications of such geometric nonabelian twisted cohomology
  of relevance in quantum field theory and in string theory is discussed in
  \cite{SSS,Schreiber}.
  \label{ReferencesOnTwisted}
\end{remark}
\begin{example}
  In particular we may consider the space of sections of the universal $\rho$-associated $V$-fiber
  $\infty$-bundle itself, hence the $\mathrm{id}_{\mathbf{B}G}$-twisted cohomology with coefficients
  in $\mathbf{c}: V/\!/G \to \mathbf{B}G$. A cocycle here is an (homotopy-)\emph{invariant}
  $$
    \xymatrix{
	  \mathbf{B}G \ar[dr]_{\mathrm{id}_{\mathbf{B}G}}^{\ }="t" \ar[rr]^{v}_{\ }="s" && V/\!/G \ar[dl]^{\mathbf{c}}
	  \\
	  & \mathbf{B}G
	  \ar@{=>} "s"; "t"
	}
  $$
  of $V$, under the $G$-action. The connected components of the hom-$\infty$-groupoid form the
  \emph{$\infty$-group cohomology} of $G$ with coefficients in $V$:
  $$
    H_{\mathrm{Grp}}(G,V)
	:=
	\pi_0 \Gamma_{\mathbf{B}G}(V/\!/G)
	\,.
  $$
  In the case where $V$ is in the image of a chain complex under the Dold-Kan correspondence,
  this statement is familiar from homological algebra: group cohomology is the derived functor
  of the invariants functor, which in turn is the hom-functor from the trivial $G$-action on the
  point (see the first item of example \ref{ExamplesOfActions} for how $\mathrm{id}_{\mathbf{B}G}$
  exhibits the trivial $G$-action on the point).
\end{example}
\begin{remark}
  More generally, of special interest is the case where
  $V$ is pointed connected, hence (by Theorem \ref{DeloopingTheorem})
  of the form $V = \mathbf{B}A$ for some $\infty$-group $A$, and so
  (by Definition \ref{cohomology}) the coefficients for degree-1 $A$-cohomology,
  and hence itself (by Theorem \ref{PrincipalInfinityBundleClassification}) the moduli $\infty$-stack for
  $A$-principal $\infty$-bundles.
  In this case $H^{[g_X]}(X, \mathbf{B}A)$ is
  \emph{degree-1 twisted $A$-cohomology}. Generally, if $V = \mathbf{B}^n A$
  it is \emph{degree-$n$ twisted $A$-cohomology}. In analogy with Definition
  \ref{cohomology} this is sometimes written
  $$
    H^{n,[g_X]}(X,A) := H^{[g_X]}(X,\mathbf{B}^n A)
	\,.
  $$

  Moreover, in this case $V/\!/G$ is itself pointed connected, hence of the form
  $\mathbf{B}\hat G$ for some $\infty$-group $\hat G$, and so the universal local
  coefficient bundle
  $$
    \xymatrix{
	  \mathbf{B}A \ar[r] & \mathbf{B}\hat G \ar[d]^{\mathbf{c}}
	  \\
	  & \mathbf{B}G
	}
  $$
  exhibits $\hat G$ as an \emph{extension of $\infty$-groups} of $G$ by $A$.
  This case we discuss below in Section \ref{ExtensionsOfCohesiveInfinityGroups}.
  \label{PointedConnectedLocalCoefficients}
\end{remark}
In this notation the local coefficient bundle $\mathbf{c}$ is left implicit.
This convenient abuse of notation is justified to some extent by the fact that there
is a \emph{universal local coefficient bundle}:
\begin{example}
  The classifying morphism of the $\mathbf{Aut}(V)$-action on some
  $V \in \mathbf{H}$ from Definition \ref{InternalAutomorphismGroup} according to
  Theorem \ref{PrincipalInfinityBundleClassification} yields a local coefficient
  $\infty$-bundle of the form
  $$
    \raisebox{20pt}{
    \xymatrix{
	  V \ar[r] & V/\!/\mathbf{Aut}(V)\ar[d]
	  \\
	  & \mathbf{B}\mathbf{Aut}(V)
	}
	}
  $$
  which we may call the \emph{universal local $V$-coefficient bundle}.
  In the case that $V$ is pointed connected and hence of the form
  $V = \mathbf{B}G$
  $$
    \raisebox{20pt}{
    \xymatrix{
	  \mathbf{B}G \ar[r] & (\mathbf{B}G)/\!/\mathbf{Aut}(\mathbf{B}G)\ar[d]
	  \\
	  & \mathbf{B}\mathbf{Aut}(\mathbf{B}G)
	}
	}
  $$
  the universal twists of the corresponding twisted $G$-cohomology are
  the \emph{$G$-$\infty$-gerbes}. These we discuss below in section
  \ref{StrucInftyGerbes}.
  \label{UniversalTwistByAutomorphisms}
\end{example}

\subsection{Extensions and twisted bundles}
\label{ExtensionsOfCohesiveInfinityGroups}

We discuss the notion of \emph{extensions} of $\infty$-groups
(see Section~\ref{StrucInftyGroups}), generalizing the traditional notion of
group extensions.
This is in fact a special case of the notion of
principal $\infty$-bundle, Definition~\ref{principalbundle}, for base space objects that
are themselves deloopings of $\infty$-groups.
For every extension of $\infty$-groups, there is the
corresponding notion of
\emph{lifts of structure $\infty$-groups} of principal $\infty$-bundles.
These are classified equivalently by trivializations of an \emph{obstruction class}
and by the twisted cohomology with coefficients in the extension itself, regarded
as a local coefficient $\infty$-bundle.

Moreover, we show
that principal $\infty$-bundles with an extended structure $\infty$-group
are equivalent to principal $\infty$-bundles with unextended structure $\infty$-group
but carrying a principal $\infty$-bundle for the \emph{extending} $\infty$-group on their
total space, which on fibers restricts to the given $\infty$-group extension.
We formalize these \emph{twisted (principal) $\infty$-bundles}
and observe that they are classified by twisted cohomology,
Definition~\ref{TwistedCohomologyInOvertopos}.

\medskip

\begin{definition}
  \label{ExtensionOfInfinityGroups}
We say a sequence of $\infty$-groups (Definition~\ref{inftygroupinootopos}),
\[
  A \to \hat G \to G
\]
in $\mathrm{Grp}(\mathbf{H})$
\emph{exhibits $\hat G$ as an extension of $G$ by $A$} if the delooping
$\mathbf{B}A \to \mathbf{B}\hat G \to \mathbf{B}G$ is a fiber sequence in $\mathbf{H}$.
\end{definition}
\begin{remark}
  By continuing this fiber sequence to the left as
  $$
    A \to \hat G \to G \to \mathbf{B}A \to \mathbf{B}\hat G \to \mathbf{B}G
  $$
  this means, by Theorem~\ref{PrincipalInfinityBundleClassification}, that
  $$
    G \simeq \hat G /\!/A
  $$
  is the quotient of the extended $\infty$-group $\hat G$ by the extending $\infty$-group $A$.
\end{remark}
\begin{definition}
  A \emph{braided $\infty$-group} is an $\infty$-group $A \in \mathrm{Grp}(\mathbf{H})$
  equipped with the following equivalent structures:
  \begin{enumerate}
    \item a lift of the defining groupal $A_\infty \simeq E_1$-action to an $E_2$-action;
	\item a group structure on the delooping $\mathbf{B}A$;
	\item a double delooping $\mathbf{B}^2 A$.
  \end{enumerate}
   \label{BraidedInfinityGroup}
\end{definition}
\begin{remark}
  The equivalence of the items in Definition \ref{BraidedInfinityGroup} is
  essentially the content of theorem 5.1.3.6 in \cite{LurieAlgebra}.
\end{remark}
\begin{definition}
For $A$ a braided $\infty$-group, Definition~\ref{BraidedInfinityGroup},
a \emph{braided-central extension} $\hat G$ of $G$ by $A$ is an extension $A \to \hat G \to G$,
Definition \ref{ExtensionOfInfinityGroups},
together with a prolongation of the defining fiber sequence one step further to the right:
$$
  \xymatrix{
    \mathbf{B}A \ar[r] &  \mathbf{B}\hat G \ar[r]^{\mathbf{p}} &  \mathbf{B}G
	\ar[r]^{\mathbf{c}} & \mathbf{B}^2 A
	\,.
  }	
$$
We write
$$
  \mathrm{Ext}(G,A) := \mathbf{H}(\mathbf{B}G, \mathbf{B}^2 A)
  \simeq (\mathbf{B}A)\mathrm{Bund}(\mathbf{B}G)
$$
for the \emph{$\infty$-groupoid of braided-central extensions} of $G$ by $A$.
\end{definition}
\begin{example}
  An ordinary group (1-group) $A$ that is braided is already abelian (by the Eckmann-Hilton argument).
  In this case a braided-central extension as above of a 1-group $G$ is a central extension of
  $G$ in the traditional sense.
\end{example}
\begin{definition}
Given an $\infty$-group extension
$\xymatrix{A \ar[r] & \hat G \ar[r]^{\Omega \mathbf{p}} & G}$
and given a $G$-principal $\infty$-bundle $P \to X$ in $\mathbf{H}$, we say that a
\emph{lift} $\hat P$ of $P$ to a $\hat G$-principal $\infty$-bundle is a lift
$\hat g_X$ of its classifying
cocycle $g_X : X \to \mathbf{B}G$, under the equivalence of Theorem
\ref{PrincipalInfinityBundleClassification}, through the extension:
$$
  \raisebox{20pt}{
  \xymatrix{
      & \mathbf{B} {\hat G}
      \ar[d]^{\mathbf{p}}
    \\
    X \ar@{-->}[ur]^{\hat g_X} \ar[r]_-{g_X} & \mathbf{B}G.
  }
  }
$$
Accordingly, the \emph{$\infty$-groupoid of lifts} of $P$ with respect to
$\mathbf{p}$ is
$$
  \mathrm{Lift}(P,\mathbf{p}) := \mathbf{H}_{/\mathbf{B}G}(g_X, \mathbf{p})
  \,.
$$
\label{PrincipalInfinityBundleExtension}
\end{definition}
\begin{observation}
  By the universal property of the $\infty$-pullback, a lift exists precisely if the
  cohomology class
  $$
    [\mathbf{c}(g_X)] := [\mathbf{c}\circ g_X] \in H^2(X, A)
  $$
  is trivial.
  \label{ObstructionClassObstructs}
\end{observation}
This is implied by Theorem \ref{ExtensionsAndTwistedCohomology}, to which we turn after
introducing the following terminology.
\begin{definition}
  In the above situation,
  we call $[\mathbf{c}(g_X)]$ the \emph{obstruction class} to the extension;
  and we call $[\mathbf{c}] \in H^2(\mathbf{B}G, A)$ the
  \emph{universal obstruction class} of extensions through $\mathbf{p}$.

  We say that a \emph{trivialization} of the obstruction cocycle
  $\mathbf{c}(g_X)$ is a homotopy $\mathbf{c}(g_X) \to *_X$ in $\mathbf{H}(X, \mathbf{B}^2 A)$
  (necessarily an equivalence),
  where ${*}_X : X \to * \to \mathbf{B}^2 A$ is the trivial cocycle. Accordingly, the
  \emph{$\infty$-groupoid of trivializations of the obstruction} is
  $$
    \mathrm{Triv}(\mathbf{c}(g_X)) := \mathbf{H}(X,\mathbf{B}^2 A)(\mathbf{c}\circ g_X, *_X)
	\,.
  $$
  \label{Obstructions}
\end{definition}
We give now three different characterizations of spaces of extensions of $\infty$-bundles.
The first two, by spaces of twisted cocycles and by spaces of trivializations
of the obstruction class, are immediate consequences of the previous discussion:
\begin{theorem}
  Let $P \to X$ be a $G$-principal $\infty$-bundle corresponding by
  Theorem \ref{PrincipalInfinityBundleClassification} to a cocycle $g_X : X \to \mathbf{B}G$.
  \begin{enumerate}
    \item
	  There is a natural equivalence
	 $$
	    \mathrm{Lift}(P, \mathbf{p}) \simeq \mathrm{Triv}(\mathbf{c}(g_X))
	 $$
	 between the $\infty$-groupoid of lifts of $P$ through $\mathbf{p}$,
     Definition \ref{PrincipalInfinityBundleExtension},	
	 and the $\infty$-groupoid of trivializations
	 of the obstruction class, Definition \ref{Obstructions}.
	 \item
	 There is a natural equivalence
	 $\mathrm{Lift}(P, \mathbf{p}) \simeq \mathbf{H}_{/\mathbf{B}G}(g_X, \mathbf{p})$
	 between the $\infty$-groupoid of lifts and the $\infty$-groupoid of
	 $g_X$-twisted cocycles relative to $\mathbf{p}$, Definition \ref{TwistedCohomologyInOvertopos},
	 hence a classification
	 $$
	   \pi_0 \mathrm{Lift}(P, \mathbf{P}) \simeq H^{1,[g_X]}(X, A)
	 $$
	 of equivalence classses of lifts by the $[g_X]$-twisted $A$-cohomology of $X$
	 relative to the local coefficient bundle
    $$
      \raisebox{20pt}{
      \xymatrix{
         \mathbf{B}A \ar[r] & \mathbf{B}\hat G \ar[d]^{\mathbf{p}}
    	 \\
	     & \mathbf{B}G\,.
       }
      }
    $$
\end{enumerate}
 \label{ExtensionsAndTwistedCohomology}
\end{theorem}
\proof
  The first statement is the special case of Lemma \ref{SliceMapsIntoPullbacks}
  where the $\infty$-pullback $E_1 \simeq f^* E_2$ in the notation there is identified
  with $\mathbf{B}\hat G \simeq \mathbf{c}^* {*}$.
  The second is evident after unwinding the definitions.
\endofproof
\begin{remark}
  For the special case that $A$ is 0-truncated, we may, by the discussion in
  \cite{NikolausWaldorf}
  identify
  $\mathbf{B}A$-principal $\infty$-bundles with $A$-\emph{bundle gerbes},
  \cite{Mur}.
  Under this identification the $\infty$-bundle classified by the
  obstruction class $[\mathbf{c}(g_X)]$ above is what is called the
  \emph{lifting bundle gerbe} of the lifting problem, see for instance
  \cite{CBMMS} for a review.
  In this case
  the first item of Theorem \ref{ExtensionsAndTwistedCohomology} reduces to
  Theorem 2.1 in \cite{Waldorf} and Theorem A (5.2.3) in \cite{NikolausWaldorf2}.
  The reduction of this statement to
  connected components,
  hence the special case of Observation \ref{ObstructionClassObstructs},
  was shown in \cite{BreenBitorseurs}.
  \label{ReferencesOnLiftings}
\end{remark}
While, therefore, the discussion of extensions of $\infty$-groups and of
lifts of structure $\infty$-groups
is just a special case of the discussion in the previous sections, this special case
admits geometric representatives of cocycles in the corresponding twisted cohomology by
twisted principal $\infty$-bundles. This we turn to now.
\begin{definition}
  \label{TwistedBundle}
  Given an extension of $\infty$-groups $A \to \hat G \xrightarrow{\Omega \mathbf{c}} G$
  and given a $G$-principal $\infty$-bundle $P \to X$, with class $[g_X] \in H^1(X,G)$,
  a \emph{$[g_X]$-twisted $A$-principal $\infty$-bundle} on $X$ is an $A$-principal
  $\infty$-bundle $\hat P \to P$ such that the cocycle $q : P \to\mathbf{B}A$
  corresponding to it under Theorem \ref{PrincipalInfinityBundleClassification}
  is a morphism of $G$-$\infty$-actions.

  The \emph{$\infty$-groupoid of $[g_X]$-twisted $A$-principal $\infty$-bundles on $X$} is
  $$
    A\mathrm{Bund}^{[g_X]}(X)
	:=
	G \mathrm{Action}( P, \mathbf{B}A )
	\subset
	\mathbf{H}(P, \mathbf{B}A)
	\,.
  $$
\end{definition}
\begin{proposition}
  \label{BundleOnTotalSpaceFromExtensionOfBundles}
  Given an $\infty$-group extension $A \to \hat G \stackrel{\Omega \mathbf{c}}{\to} G$,
  an extension of a $G$-principal $\infty$-bundle $P \to X$
  to a $\hat G$-principal $\infty$-bundle, Definition \ref{PrincipalInfinityBundleExtension}, induces
  an $A$-principal $\infty$-bundle $\hat P \to P$
  fitting into a pasting diagram of $\infty$-pullbacks of the form
$$
  \raisebox{30pt}{
  \xymatrix{
    \hat G \ar[r] \ar[d]^{\Omega \mathbf{c}} & \hat P \ar[r] \ar[d] &  {*} \ar[d]
    \\
    G \ar[r] \ar[d] & P \ar[r]^{q} \ar[d]& \mathbf{B}A \ar[r] \ar[d] & {*} \ar[d]
    \\
    {*} \ar[r]^{x}&
    X \ar@/_1pc/[rr]_g \ar[r]^{\hat g} & \mathbf{B}\hat G
    \ar[r]^{\mathbf{c}}& \mathbf{B}G.
  }
  }
$$
  In particular, it has the following properties:
  \begin{enumerate}
    \item $\hat P \to P$ is a $[g_X]$-twisted $A$-principal bundle, Definition \ref{TwistedBundle};
	\item for all points $x : * \to X$ the restriction of $\hat P \to P$ to the
	fiber $P_x$ is canonically equivalent to the $\infty$-group extension $\hat G \to G$.
  \end{enumerate}
\end{proposition}
\proof
This follows from repeated application of the pasting law for $\infty$-pullbacks,
Proposition~\ref{PastingLawForPullbacks}.

The bottom composite $g : X \to \mathbf{B}G$ is a cocycle for the given $G$-principal $\infty$-bundle
$P \to X$ and it factors through $\hat g : X \to \mathbf{B}\hat G$ by assumption of the existence of
the extension $\hat P \to P$.

Since also the bottom right square is an $\infty$-pullback by the given $\infty$-group extension,
the pasting law asserts that the square over $\hat g$ is also an  $\infty$-pullback, and then that so
is the square over $q$. This exhibits $\hat P$ as an $A$-principal $\infty$-bundle over $P$ classified
by the cocycle $q$ on $P$. By Corollary~\ref{SectionsOfAssociatedIsActionMorphism}
this $\hat P \to P$ is twisted $G$-equivariant.

Now choose any point $x : {*} \to X$ of the base space
as on the left of the diagram. Pulling this back upwards
through the diagram and using the pasting law and the
definition of loop space objects $G \simeq \Omega
\mathbf{B}G \simeq * \times_{\mathbf{B}G} *$ the diagram
completes by $\infty$-pullback squares on the left as indicated, which proves the claim.
\endofproof

\begin{theorem}
  The construction of Proposition~\ref{BundleOnTotalSpaceFromExtensionOfBundles}
  extends to an equivalence of $\infty$-groupoids
  $$
    A\mathrm{Bund}^{[g_X]}(X)
	\simeq
	\mathbf{H}_{/\mathbf{B}G}(g_X, \mathbf{c})
  $$
  between that of $[g_X]$-twisted $A$-principal bundles on $X$,
  Definition \ref{TwistedBundle}, and the cocycle $\infty$-groupoid of
  degree-1 $[g_X]$-twisted $A$-cohomology, Definition \ref{TwistedCohomologyInOvertopos}.

  In particular the classification of $[g_X]$-twisted $A$-principal bundles is
  $$
    A\mathrm{Bund}^{[g_X]}(X)_{/\sim}
	\simeq
	H^{1,[g_X]}(X, A)
	\,.
  $$
  \label{ClassificationOfTwistedGEquivariantBundles}
\end{theorem}
\proof
  For $G = *$ the trivial group, the statement reduces to
  Theorem~\ref{PrincipalInfinityBundleClassification}. The general proof
  works along the same lines as the proof of that theorem.
  The key step is the generalization of the proof of Proposition~\ref{LocalTrivialityImpliesCocycle}.
  This proceeds verbatim as there, only with $\mathrm{pt} : * \to \mathbf{B}G$ generalized to
  $i : \mathbf{B}A \to \mathbf{B}\hat G$.
  The morphism of $G$-actions
  $P \to \mathbf{B}A$ and a choice of effective epimorphism $U \to X$ over
  which $P \to X$ trivializes gives rise to a morphism in
  $\mathbf{H}^{\Delta[1]}_{/(* \to \mathbf{B}G)}$ which involves the diagram
    $$
    \raisebox{20pt}{
    \xymatrix{
	  U \times G \ar@{->>}[r] \ar[d] & P \ar[r] \ar[d] & \mathbf{B}A \ar[d]^{i}
	  \\
	  U \ar@{->>}[r] & X \ar[r] & \mathbf{B}\hat G
	}
	}
	\;\;
	\simeq
	\;\;
    \raisebox{20pt}{
    \xymatrix{
	  U \times G \ar@{->>}[rr] \ar[d] & & \mathbf{B}A \ar[d]^{i}
	  \\
	  U \ar[r] & {*} \ar[r]^{\mathrm{pt}} & \mathbf{B}\hat G
	}
	}
  $$
  in $\mathbf{H}$. (We are using that for the 0-connected object $\mathbf{B}\hat G$
  every morphism $* \to \mathbf{B}G$ factors through $\mathbf{B}\hat G \to \mathbf{B}G$.)
  Here the total rectangle and the left square on the left are $\infty$-pullbacks,
  and we need to show that the right square on the left is then also an
  $\infty$-pullback. Notice that by the pasting law the rectangle on the
  right is indeed equivalent to the pasting of $\infty$-pullbacks
  $$
    \raisebox{20pt}{
    \xymatrix{
	  U \times G \ar@{->}[r] \ar[d] & G \ar[r] \ar[d] &  \mathbf{B}A \ar[d]^{i}
	  \\
	  U \ar[r] & {*} \ar[r]^{\mathrm{pt}} & \mathbf{B}\hat G
	}
	}
  $$
  so that the relation
  $$
    U^{\times^{n+1}_X} \times G
	\simeq
	i^* (U^{\times^{n+1}_X})
  $$
  holds. With this the proof finishes as in the proof of
  Proposition~\ref{LocalTrivialityImpliesCocycle}, with $\mathrm{pt}^*$
  generalized to $i^*$.
\endofproof

\begin{remark}
  Aspects of special cases of this theorem can be identified in the literature.
  For the special case of ordinary extensions of ordinary Lie groups,
  the equivalence of the corresponding extensions of a principal bundle with
  certain equivariant structures on its total space is essentially the content
  of \cite{Mackenzie, Androulidakis}.
  In particular the  twisted unitary bundles or \emph{gerbe modules}
  of twisted K-theory \cite{CBMMS}
  are equivalent to such structures.

  For the case of $\mathbf{B}U(1)$-extensions of Lie groups, such as the
  $\mathrm{String}$-2-group, the equivalence of the corresponding
  $\mathrm{String}$-principal 2-bundles,
  by the above theorem, to certain bundle gerbes on the total spaces of
  principal bundles underlies constructions such as in \cite{Redden}.
  Similarly the bundle gerbes on double covers considered in
  \cite{SSW} are $\mathbf{B}U(1)$-principal 2-bundles on $\mathbb{Z}_2$-principal
  bundles arising by the above theorem from the extension
  $\mathbf{B}U(1) \to \mathbf{Aut}(\mathbf{B}U(1)) \to \mathbb{Z}_2$,
  a special case of the extensions that we consider in the next Section
  \ref{StrucInftyGerbes}.
  \label{ReferencesTwistedBundles}
\end{remark}

\subsection{Gerbes}
\label{section.InfinityGerbes}
\label{StrucInftyGerbes}

Recall from Remark \ref{PointedConnectedLocalCoefficients} above
that in an $\infty$-topos $\mathbf{H}$, those $V$-fiber $\infty$-bundles
$E\to X$ whose typical fiber $V$ is pointed connected are of
special relevance.  Recall that such a $V$ is the moduli $\infty$-stack $V = \mathbf{B}G$
of $G$-principal $\infty$-bundles for some $\infty$-group $G$.
Due to their local triviality,
when regarded as objects in the slice $\infty$-topos $\mathbf{H}_{/X}$,
these $\mathbf{B}G$-fiber $\infty$-bundles are themselves \emph{connected objects}.
Generally, for $\mathcal{X}$ an $\infty$-topos regarded as an $\infty$-topos
of $\infty$-stacks over a given space $X$, it makes sense to consider
its connected objects as $\infty$-bundles over $X$. Here we discuss these
\emph{$\infty$-gerbes}.

\medskip

In the following discussion it is useful to consider two $\infty$-toposes:
\begin{enumerate}
  \item an ``ambient'' $\infty$-topos $\mathbf{H}$ as before, to be thought of
    as an $\infty$-topos ``of all geometric homotopy types'' for a given notion of geometry
    (recall the discussion in section \ref{Overview}),
	in which $\infty$-bundles are given by \emph{morphisms} and the terminal
	object plays the role of the geometric point $*$;
  \item an $\infty$-topos $\mathcal{X}$, to be thought of as the topos-theoretic
   incarnation of a single geometric homotopy type (space) $X$, hence as
   an $\infty$-topos of
  ``geometric homotopy types {\'e}tale over $X$'', in which an $\infty$-bundle over $X$
  is given by an \emph{object} and the terminal object plays the role of the base space $X$.

  In practice, $\mathcal{X}$ is the slice
  $\mathbf{H}_{/X}$ of the previous  ambient $\infty$-topos over $X \in \mathbf{H}$, or
  the smaller $\infty$-topos $\mathcal{X} = \mathrm{Sh}_\infty(X)$ of (internal)
  $\infty$-stacks over $X$ (hence {\'e}tale objects over $X$, see section 3.10.7 of \cite{Schreiber}).
\end{enumerate}
In topos-theory literature the role of $\mathbf{H}$ above is sometimes referred to as that of a
\emph{gros} topos and then the role of $\mathcal{X}$ is referred to as that of a \emph{petit} topos.
The reader should beware that much of the classical literature on gerbes is written
from the point of view of only the \emph{petit} topos $\mathcal{X}$.
For the following, recall remark \ref{CohomologyOverX} on cohomology in slice toposes.

\medskip

The original definition of a \emph{gerbe}
on $X$ as given by \cite{Giraud} is: a stack $E$ (i.e.\  a
1-truncated $\infty$-stack) over $X$ that is
1. \emph{locally non-empty} and 2. \emph{locally connected}.
In the more intrinsic language of higher topos theory, these
two conditions simply say that $E$ is
a \emph{connected object} (Definition 6.5.1.10 in \cite{Lurie}):
1. the terminal morphism $E \to *$ is an
effective epimorphism and 2. the 0th homotopy sheaf is trivial,
$\pi_0(E) \simeq *$. This reformulation is
made explicit in the literature for instance in
Section 5 of \cite{JardineLuo} and in Section 7.2.2 of \cite{Lurie}. Therefore:
\begin{definition}
  For $\mathcal{X}$ an $\infty$-topos, a \emph{gerbe} in $\mathcal{X}$ is an object $E \in \mathcal{X}$ which is
  \begin{enumerate}
    \item connected;
	\item 1-truncated.
  \end{enumerate}
  For $X \in \mathbf{H}$ an object, a \emph{gerbe $E$ over $X$} is a gerbe in the slice
  $\mathbf{H}_{/X}$. This is an object $E \in \mathbf{H}$ together with an effective epimorphism $E \to X$ such that $\pi_i(E) = X$ for all $i \neq 1$.
  \label{Gerbe}
\end{definition}
\begin{remark}
  Notice that conceptually this is different from the
  notion of \emph{bundle gerbe} introduced in \cite{Mur}
  (see \cite{NikolausWaldorf} for a review).
  Bundle gerbes are presentations of \emph{principal} $\infty$-bundles
  (Definition \ref{principalbundle}).
  But gerbes -- at least the \emph{$G$-gerbes}
  considered in a moment in Definition \ref{G-Gerbe} --
  are $V$-fiber $\infty$-bundles (Definition \ref{FiberBundle})
  hence \emph{associated} to principal $\infty$-bundles
  (Proposition \ref{VBundleIsAssociated}) with the special
  property of having pointed connected fibers.
  By Theorem \ref{VBundleClassification}
  $V$-fiber $\infty$-bundles may be identified with their underlying $\mathbf{Aut}(V)$-principal
  $\infty$-bundles and so one may identify $G$-gerbes with nonabelian
  $\mathrm{Aut}(\mathbf{B}G)$-bundle gerbes
  (see also around Corollary~\ref{ClassificationOfGGerbes} below),
  but considered generally, neither of these
  two notions is a special case of the other. Therefore the terminology is
  slightly unfortunate, but it is standard.
\end{remark}
Definition \ref{Gerbe} has various obvious generalizations. The following is considered in \cite{Lurie}.
\begin{definition}
  \index{$\infty$-gerbe!EM $n$-gerbe}
 For $n \in \mathbb{N}$, an \emph{EM $n$-gerbe} is an object $E \in \mathcal{X}$ which is
 \begin{enumerate}
   \item $(n-1)$-connected;
   \item $n$-truncated.
 \end{enumerate}
\end{definition}
\begin{remark}
This is almost the definition of an
\emph{Eilenberg-Mac Lane object} in $\mathcal{X}$,
only that the condition requiring a global section
$* \to E$ (hence $X \to E$) is missing. Indeed, the
Eilenberg-Mac Lane objects of degree $n$ in
$\mathcal{X}$ are precisely the EM $n$-gerbes of
\emph{trivial class}, according to Corollary~\ref{ClassificationOfGGerbes} below.
\end{remark}
There is also an earlier established definition of
\emph{2-gerbes} in the literature \cite{Breen}, which is
more general than EM 2-gerbes. Stated in the above fashion it reads as follows.
\begin{definition}[Breen \cite{Breen}]
  \index{$\infty$-gerbe!2-gerbe}
  \index{2-gerbe}
  A \emph{2-gerbe} in $\mathcal{X}$ is an object $E \in \mathcal{X}$ which is
  \begin{enumerate}
    \item connected;
	\item 2-truncated.
  \end{enumerate}
\end{definition}
This definition has an evident generalization to arbitrary degree, which we adopt here.
\begin{definition}
  \label{nGerbe}
  An \emph{$n$-gerbe} in $\mathcal{X}$ is an object $E \in \mathcal{X}$ which is
   \begin{enumerate}
      \item connected;
	  \item $n$-truncated.
   \end{enumerate}
In particular an \emph{$\infty$-gerbe} is a connected object.
\end{definition}
The real interest is in those $\infty$-gerbes which have a prescribed
typical fiber:
\begin{remark}
By the above, $\infty$-gerbes (and hence EM $n$-gerbes
and 2-gerbes and hence gerbes) are much like deloopings
of $\infty$-groups (Theorem \ref{DeloopingTheorem})
only that there is no requirement that there
exists a global section. An $\infty$-gerbe for which there exists a
global section $X \to E$ is called \emph{trivializable}. By
Theorem~\ref{DeloopingTheorem} trivializable $\infty$-gerbes are equivalent to $\infty$-group objects in $\mathcal{X}$
(and the $\infty$-groupoids of all of these are equivalent
when transformations are required to preserve the canonical global section).
\end{remark}
But \emph{stalkwise} every $\infty$-gerbe $E$ is of this form. For let
$$
  (x^* \dashv x_*) :
  \xymatrix{
     \mathrm{Grpd}_{\infty}
	    \ar@{<-}@<+3pt>[r]^-{x^*}
	    \ar@<-3pt>[r]_-{x_*}
		&
     \mathcal{X}	
  }
$$
be a topos point. Then the stalk $x^* E \in \mathrm{Grpd}_{\infty}$
of the $\infty$-gerbe is connected: because inverse images
preserve the finite $\infty$-limits involved in the definition of homotopy sheaves, and preserve the terminal object. Therefore
$$
  \pi_0\,  x^* E \simeq x^* \pi_0 E \simeq x^* * \simeq *
  \,.
$$
Hence for every point $x$ we have a stalk $\infty$-group $G_x$ and an equivalence
$$
  x^* E \simeq B G_x
  \,.
$$
Therefore one is interested in the following notion.
\begin{definition}
  For $G \in \mathrm{Grp}(\mathcal{X})$ an $\infty$-group
  object, a \emph{$G$-$\infty$-gerbe} is an $\infty$-gerbe $E$ such that there exists
  \begin{enumerate}
    \item an effective epimorphism $\xymatrix{U \ar@{->>}[r] & X}$ (onto the terminal object $X$ of $\mathcal{X}$);
	\item an equivalence $E|_U \simeq \mathbf{B} G|_U$.
  \end{enumerate}
  Equivalently: a $G$-$\infty$-gerbe is a $\mathbf{B}G$-fiber $\infty$-bundle over the terminal object $X$ of $\mathcal{X}$,
  according to Definition \ref{FiberBundle}.
  \label{G-Gerbe}
\end{definition}
In words this says that a $G$-$\infty$-gerbe is one that locally looks like the
moduli $\infty$-stack of $G$-principal $\infty$-bundles.
\begin{example}
  For $X$ a topological space and $\mathcal{X} =
  \mathrm{Sh}_\infty(X)$ the $\infty$-topos of $\infty$-sheaves over it, these notions reduce to the following.
  \begin{itemize}
    \item a 0-group object $G \in \tau_{0}\mathrm{Grp}(\mathcal{X})
    \subset \mathrm{Grp}(\mathcal{X})$ is a sheaf of groups on $X$
    (here $\tau_0\mathrm{Grp}(\mathcal{X})$ denotes the 0-truncation of
    $\mathrm{Grp}(\mathcal{X})$);
	\item for $\{U_i \to X\}$ any open cover, the canonical
	morphism $\coprod_i U_i \to X$ is an effective epimorphism to the terminal object;
	\item $(\mathbf{B}G)|_{U_i}$ is the stack of
	   $G|_{U_i}$-principal bundles ($G|_{U_i}$-torsors).
  \end{itemize}
\end{example}
It is clear that one way to construct a $G$-$\infty$-gerbe
should be to start with an $\mathbf{Aut}(\mathbf{B}G)$-principal
$\infty$-bundle, Remark \ref{UniversalTwistByAutomorphisms},
and then canonically \emph{associate} a fiber $\infty$-bundle to it.
\begin{example}
  \label{automorphism2GroupAbstractly}
  For $G \in \tau_{0}\mathrm{Grp}(\mathrm{Grpd}_{\infty})$ an
  ordinary group, $\mathbf{Aut}(\mathbf{B}G)$ is usually called the
  \emph{automorphism 2-group} of $G$. Its underlying groupoid is equivalent to
  \[
  \mathbf{Aut}(G) \times G\rightrightarrows
  \mathbf{Aut}(G),
  \]
  the action groupoid for the action of $G$ on $\mathbf{Aut}(G)$
  via the homomorphism $\mathrm{Ad}\colon G\to \mathbf{Aut}(G)$.
\end{example}
\begin{corollary}
  Let $\mathcal{X}$ be an $\infty$-topos.
  Then for $G \in \mathrm{Grp}(\mathcal{X})$ any
  $\infty$-group object, $G$-$\infty$-gerbes are classified by $\mathbf{Aut}(\mathbf{B}G)$-cohomology:
  $$
    \pi_0 G \mathrm{Gerbe}
	\simeq
	\pi_0 \mathcal{X}(X,\mathbf{B}\mathbf{Aut}(\mathbf{B}G))
	=:
	H^1(X,\mathbf{Aut}(\mathbf{B}G))
	\,.
  $$
  \label{ClassificationOfGGerbes}
\end{corollary}
\proof
  This is the special case of Theorem \ref{VBundleClassification}
  for $V = \mathbf{B}G$.
\endofproof
For the case that $G$ is 0-truncated (an ordinary group object)
this is the content of Theorem 23 in \cite{JardineLuo}.
\begin{example}
  For $G \in \tau_{\leq 0}\mathrm{Grp}(\mathcal{X}) \subset
  \mathrm{Grp}(\mathcal{X})$ an ordinary 1-group object,
  this reproduces the classical result of \cite{Giraud},
  which originally motivated the whole subject:
  by Example~\ref{automorphism2GroupAbstractly} in this case $\mathbf{Aut}(\mathbf{B}G)$ is the traditional automorphism 2-group and
  $H^1(X, \mathbf{Aut}(\mathbf{B}G))$
  is Giraud's nonabelian $G$-cohomology that classifies $G$-gerbes
  (for arbitrary \emph{band}, see Definition \ref{BandOfInfinityGerbe} below).

  For $G \in \tau_{\leq 1}\mathrm{Grp}(\mathcal{X}) \subset
  \mathrm{Grp}(\mathcal{X})$ a 2-group, we recover the classification of 2-gerbes
  as in \cite{Breen,BreenNotes}.
\end{example}
\begin{remark}
  In Section 7.2.2 of \cite{Lurie} the special case that here
  we called \emph{EM-$n$-gerbes} is considered. Beware
  that there are further differences: for instance the notion
  of morphisms between $n$-gerbes as defined in \cite{Lurie} is more
  restrictive than the notion considered here. For instance with our definition
  (and hence also that in \cite{Breen}) each group automorphism of
  an abelian group object $A$ induces an automorphism of the
  trivial $A$-2-gerbe $\mathbf{B}^2 A$. But, except for the identity,
  this is not admitted in \cite{Lurie} (manifestly so by the diagram
  above Lemma 7.2.2.24 there). Accordingly, the classification
  result in \cite{Lurie} is different: it involves  the cohomology group
  $H^{n+1}(X, A)$. Notice that there is a canonical morphism
  $$
    H^{n+1}(X, A) \to H^1(X, \mathbf{Aut}(\mathbf{B}^n A))
  $$
  induced from the morphism $\mathbf{B}^{n+1}A \to \mathbf{Aut}(\mathbf{B}^n A)$.
\end{remark}

We now discuss how the $\infty$-group extensions (Definition \ref{ExtensionOfInfinityGroups})
given by the Postnikov stages of $\mathbf{Aut}(\mathbf{B}G)$, induce the notion of
\emph{band} of a gerbe, and how the corresponding twisted cohomology,
according to Remark \ref{ExtensionsAndTwistedCohomology}, reproduces
the original definition of nonabelian cohomology in \cite{Giraud} and generalizes
it to higher degree.
\begin{definition}
 \label{outerAutomorphismInfinityGroup}
 Fix $k \in \mathbb{N}$. For $G \in \infty\mathrm{Grp}(\mathcal{X})$ a
 $k$-truncated $\infty$-group object (a $(k+1)$-group), write
 $$
   \mathbf{Out}(G) := \tau_{k}\mathbf{Aut}(\mathbf{B}G)
 $$
 for the $k$-truncation of $\mathbf{Aut}(\mathbf{B}G)$. (Notice
 that this is still an $\infty$-group, since by Lemma 6.5.1.2 in
 \cite{Lurie} $\tau_n$ preserves all $\infty$-colimits and additionally
 all products.) We call this the \emph{outer automorphism $n$-group} of $G$.

 In other words, we write
 $$
   \mathbf{c} : \mathbf{B}\mathbf{Aut}(\mathbf{B}G) \to \mathbf{B}\mathbf{Out}(G)
 $$
 for the top Postnikov stage of $\mathbf{B}\mathbf{Aut}(\mathbf{B}G)$.
\end{definition}
\begin{example}
  Let $G \in \tau_0\mathrm{Grp}(\mathrm{Grpd}_{\infty})$ be
  a 0-truncated group object, an
  ordinary group.
  Then by Example \ref{automorphism2GroupAbstractly},
  $\mathbf{Out}(G)$ is the coimage of
  $\mathrm{Ad} : G \to \mathrm{Aut}(G)$, which is the traditional group of outer automorphisms of $G$.
\end{example}
\begin{definition}
  Write $\mathbf{B}^2 \mathbf{Z}(G)$ for the $\infty$-fiber of
  the morphism $\mathbf{c}$ from Definition \ref{outerAutomorphismInfinityGroup},
  fitting into a fiber sequence
  $$
    \xymatrix{
      \mathbf{B}^2 \mathbf{Z}(G) \ar[r] &
	    \mathbf{B}\mathbf{Aut}(\mathbf{B}G) \ar[d]^{\mathbf{c}}
	  \\
	  & \mathbf{B}\mathbf{Out}(G)
	}
	\,.
  $$
  We call $\mathbf{Z}(G)$ the \emph{braided center} of the $\infty$-group $G$.
  \label{Center}
\end{definition}
\begin{remark}
  To see that the fiber of $\Omega \mathbf{c}$ here is indeed the delooping of
  a group, notice that by theorem \ref{DeloopingTheorem} one has to see that
  it is connected and pointed.
  Now the fiber of $\Omega \mathbf{c}$ is connected due to
  definition of $\mathbf{c}$ as a truncation map and the induced
  long exact sequence of (sheaves of) homotopy groups. It is moreover
  pointed since $\Omega \mathbf{c}$, being a morphism    of groups, is a pointed morphism
  (the point being the neutral element) and using the universal property of the
  homotopy fiber.
\end{remark}
\begin{example}
  For $G$ an ordinary group, so that $\mathbf{Aut}(\mathbf{B}G)$
  is the automorphism 2-group from Example \ref{automorphism2GroupAbstractly},
  $\mathbf{Z}(G)$ is the center of $G$ in the traditional sense.
\end{example}
By Corollary \ref{ClassificationOfGGerbes} there is an  induced morphism
$$
  \mathrm{Band}
    :
  \pi_0 G \mathrm{Gerbe} \to H^1(X, \mathbf{Out}(G))
   \,.
$$
\begin{definition}
  \label{BandOfInfinityGerbe}
  For $E \in G \mathrm{Gerbe}$ we call $\mathrm{Band}(E)$ the \emph{band} of $E$.

  By using Definition \ref{Center} in Definition \ref{TwistedCohomologyInOvertopos},
  given a band $[\phi_X] \in H^1(X, \mathbf{Out}(G))$,
  we may regard it as a twist for twisted $\mathbf{Z}(G)$-cohomology,
  classifying $G$-gerbes with this band:
  $$
    \pi_0 G \mathrm{Gerbe}^{[\phi_X]}(X) \simeq H^{2,[\phi_X]}(X, \mathbf{Z}(G))
	\,.
  $$
\end{definition}
\begin{remark}
  The original definition of \emph{gerbe with band} in \cite{Giraud} is slightly
  more general than that of \emph{$G$-gerbe} (with band) in \cite{Breen}: in the former
  the local sheaf of groups whose delooping is locally equivalent to the gerbe need not
  descend to the base. These more general Giraud gerbes are 1-gerbes in the sense of
  Definition \ref{nGerbe}, but only the slightly more restrictive $G$-gerbes of Breen
  have the good property of being connected fiber $\infty$-bundles. From our
  perspective this is the decisive property of gerbes, and the notion of band is
  relevant only in this case.
\end{remark}
\begin{example}
  For $G$ a 0-group this reduces to the notion of band as introduced in  \cite{Giraud},
  for the case of $G$-gerbes as in \cite{Breen}.
\end{example}

\medskip

\noindent{\bf Acknowledgements.}
The writeup of this article and the companion \cite{NSSb} was
initiated during a visit by the first two authors to the third author's institution,
University of Glasgow, in summer 2011.  It was completed in summer 2012
when all three authors were guests at the
Erwin Schr\"{o}dinger Institute in Vienna.
The authors gratefully acknowledge the support of
the Engineering and Physical Sciences Research Council
grant number EP/I010610/1 and the support of the ESI;
D.S. gratefully acknowledges the support of the Australian Research Council
(grant number DP120100106); U.S. acknowledges the support of
the Dutch Research Organization
NWO (project number 613.000.802).
U.S. thanks Domenico Fiorenza for inspiring discussion about
twisted cohomology. We also thank Eivind O. Hjelle for pointing out a missing assumption in the statement of Lemma \ref{TrivialityOfGBundlesOverPoints} in an earlier version.

\addcontentsline{toc}{section}{References}

\end{document}